
\catcode'32=9
\magnification=1200

\voffset=1cm

\font\tenpc=cmcsc10

\font\eightrm=cmr8
\font\eighti=cmmi8
\font\eightsy=cmsy8
\font\eightbf=cmbx8
\font\eighttt=cmtt8
\font\eightit=cmti8
\font\eightsl=cmsl8
\font\sixrm=cmr6
\font\sixi=cmmi6
\font\sixsy=cmsy6
\font\sixbf=cmbx6

\skewchar\eighti='177 \skewchar\sixi='177
\skewchar\eightsy='60 \skewchar\sixsy='60

\font\tengoth=eufm10
\font\tenbboard=msbm10
\font\eightgoth=eufm7 at 8pt
\font\eightbboard=msbm7 at 8pt
\font\sevengoth=eufm7
\font\sevenbboard=msbm7
\font\sixgoth=eufm5 at 6 pt
\font\fivegoth=eufm5

\font\tengoth=eufm10
\font\tenbboard=msbm10
\font\eightgoth=eufm7 at 8pt
\font\eightbboard=msbm7 at 8pt
\font\sevengoth=eufm7
\font\sevenbboard=msbm7
\font\sixgoth=eufm5 at 6 pt
\font\fivegoth=eufm5

\newfam\gothfam
\newfam\bboardfam

\catcode`\@=11

\def\raggedbottom{\topskip 10pt plus 36pt
\r@ggedbottomtrue}
\def\pc#1#2|{{\bigf@ntpc #1\penalty
\@MM\hskip\z@skip\smallf@ntpc #2}}

\def\tenpoint{%
  \textfont0=\tenrm \scriptfont0=\sevenrm \scriptscriptfont0=\fiverm
  \def\rm{\fam\z@\tenrm}%
  \textfont1=\teni \scriptfont1=\seveni \scriptscriptfont1=\fivei
  \def\oldstyle{\fam\@ne\teni}%
  \textfont2=\tensy \scriptfont2=\sevensy \scriptscriptfont2=\fivesy
  \textfont\gothfam=\tengoth \scriptfont\gothfam=\sevengoth
  \scriptscriptfont\gothfam=\fivegoth
  \def\goth{\fam\gothfam\tengoth}%
  \textfont\bboardfam=\tenbboard \scriptfont\bboardfam=\sevenbboard
  \scriptscriptfont\bboardfam=\sevenbboard
  \def\bboard{\fam\bboardfam}%
  \textfont\itfam=\tenit
  \def\it{\fam\itfam\tenit}%
  \textfont\slfam=\tensl
  \def\sl{\fam\slfam\tensl}%
  \textfont\bffam=\tenbf \scriptfont\bffam=\sevenbf
  \scriptscriptfont\bffam=\fivebf
  \def\bf{\fam\bffam\tenbf}%
  \textfont\ttfam=\tentt
  \def\tt{\fam\ttfam\tentt}%
  \abovedisplayskip=12pt plus 3pt minus 9pt
  \abovedisplayshortskip=0pt plus 3pt
  \belowdisplayskip=12pt plus 3pt minus 9pt
  \belowdisplayshortskip=7pt plus 3pt minus 4pt
  \smallskipamount=3pt plus 1pt minus 1pt
  \medskipamount=6pt plus 2pt minus 2pt
  \bigskipamount=12pt plus 4pt minus 4pt
  \normalbaselineskip=12pt
  \setbox\strutbox=\hbox{\vrule height8.5pt depth3.5pt width0pt}%
  \let\bigf@ntpc=\tenrm \let\smallf@ntpc=\sevenrm
  \let\petcap=\tenpc
  \normalbaselines\rm}
\def\eightpoint{%
  \textfont0=\eightrm \scriptfont0=\sixrm \scriptscriptfont0=\fiverm
  \def\rm{\fam\z@\eightrm}%
  \textfont1=\eighti \scriptfont1=\sixi \scriptscriptfont1=\fivei
  \def\oldstyle{\fam\@ne\eighti}%
  \textfont2=\eightsy \scriptfont2=\sixsy \scriptscriptfont2=\fivesy
  \textfont\gothfam=\eightgoth \scriptfont\gothfam=\sixgoth
  \scriptscriptfont\gothfam=\fivegoth
  \def\goth{\fam\gothfam\eightgoth}%
  \textfont\bboardfam=\eightbboard \scriptfont\bboardfam=\sevenbboard
  \scriptscriptfont\bboardfam=\sevenbboard
  \def\bboard{\fam\bboardfam}%
  \textfont\itfam=\eightit
  \def\it{\fam\itfam\eightit}%
  \textfont\slfam=\eightsl
  \def\sl{\fam\slfam\eightsl}%
  \textfont\bffam=\eightbf \scriptfont\bffam=\sixbf
  \scriptscriptfont\bffam=\fivebf
  \def\bf{\fam\bffam\eightbf}%
  \textfont\ttfam=\eighttt
  \def\tt{\fam\ttfam\eighttt}%
  \abovedisplayskip=9pt plus 2pt minus 6pt
  \abovedisplayshortskip=0pt plus 2pt
  \belowdisplayskip=9pt plus 2pt minus 6pt
  \belowdisplayshortskip=5pt plus 2pt minus 3pt
  \smallskipamount=2pt plus 1pt minus 1pt
  \medskipamount=4pt plus 2pt minus 1pt
  \bigskipamount=9pt plus 3pt minus 3pt
  \normalbaselineskip=9pt
  \setbox\strutbox=\hbox{\vrule height7pt depth2pt width0pt}%
  \let\bigf@ntpc=\eightrm \let\smallf@ntpc=\sixrm
  \normalbaselines\rm}

\tenpoint

\frenchspacing


\newif\ifpagetitre
\newtoks\auteurcourant \auteurcourant={\hfil}
\newtoks\titrecourant \titrecourant={\hfil}

\def\appeln@te{}
\def\vfootnote#1{\def\@parameter{#1}\insert\footins\bgroup\eightpoint
  \interlinepenalty\interfootnotelinepenalty
  \splittopskip\ht\strutbox 
  \splitmaxdepth\dp\strutbox \floatingpenalty\@MM
  \leftskip\z@skip \rightskip\z@skip
  \ifx\appeln@te\@parameter\indent \else{\noindent #1\ }\fi
  \footstrut\futurelet\next\fo@t}

\pretolerance=500 \tolerance=1000 \brokenpenalty=5000
\newdimen\hmargehaute \hmargehaute=0cm
\newdimen\lpage \lpage=13.3cm
\newdimen\hpage \hpage=20cm
\newdimen\lmargeext \lmargeext=1cm
\hsize=11.25cm
\vsize=18cm
\parskip 0pt
\parindent=12pt

\def\margehaute{\vbox to \hmargehaute{\vss}}%
\def\margebasse{\vss}

\output{\shipout\vbox to \hpage{\margehaute\nointerlineskip
  \corpsdepage\margebasse}
  \advancepageno \global\pagetitrefalse
  \ifnum\outputpenalty>-20000 \else\dosupereject\fi}

\def\corpsdepage{\hbox to \lpage{\hss\pagetexte\hskip\lmargeext}}
\def\pagetexte{\vbox{\makeheadline\pagebody\makefootline}}
\headline={\ifpagetitre\titleheadline \else
  \ifodd\pageno\rightheadline \else\leftheadline\fi\fi}
\def\leftheadline{\eightpoint\hfil\the\auteurcourant\hfil}
\def\rightheadline{\eightpoint\hfil\the\titrecourant\hfil}
\def\titleheadline{\hfill}
\pagetitretrue

\def\footnoterule{\kern-6\p@
  \hrule width 2truein \kern 5.6\p@} 

\def\pd#1#2 {\pc#1#2| }

\def\pointir{\discretionary{.}{}{.\kern.35em---\kern.7em}\nobreak
\hskip 0em plus .3em minus .4em }

\def\abstract#1{\vbox{\eightpoint \pc ABSTRACT|\pointir #1}}

\def\titre#1|{\message{#1}
              \par\vskip 30pt plus 24pt minus 3pt\penalty -1000
              \vskip 0pt plus -24pt minus 3pt\penalty -1000
              \centerline{\bf #1}
              \vskip 5pt
              \penalty 10000 }

\def\section#1|{\par\vskip .3cm
                {\bf #1}\pointir}

\def\ssection#1|{\par\vskip .2cm
                {\it #1}\pointir}

\long\def\th#1|#2\finth{\par\medskip
              {\petcap #1\pointir}{\it #2}\par\smallskip}

\long\def\tha#1|#2\fintha{\par\medskip
                    {\petcap #1.}\par\nobreak{\it #2}\par\smallskip}

\def\rem#1|{\par\medskip
            {{\it #1}.\quad}}

\def\rema#1|{\par\medskip
             {{\it #1.}\par\nobreak }}

\def\article#1|#2|#3|#4|#5|#6|#7|
    {{\leftskip=7mm\noindent
     \hangindent=2mm\hangafter=1
     \llap{[#1]\hskip.35em}{#2}.\quad
     #3, {\sl #4}, vol.\nobreak\ {\bf #5}, {\oldstyle #6},
     p.\nobreak\ #7.\par}}
\def\livre#1|#2|#3|#4|
    {{\leftskip=7mm\noindent
    \hangindent=2mm\hangafter=1
    \llap{[#1]\hskip.35em}{#2}.\quad
    {\sl #3}.\quad #4.\par}}
\def\divers#1|#2|#3|
    {{\leftskip=7mm\noindent
    \hangindent=2mm\hangafter=1
     \llap{[#1]\hskip.35em}{#2}.\quad
     #3.\par}}
\mathchardef\conj="0365
\def\proof{\par{\it Proof}.\quad}
\def\qed{\quad\raise -2pt\hbox{\vrule\vbox to 10pt{\hrule width 4pt
\vfill\hrule}\vrule}}

\def\cqfd{\penalty 500 \hbox{\qed}\par\smallskip}
\def\decale#1|{\par\noindent\hskip 28pt\llap{#1}\kern 5pt}

\catcode`\@=12


\catcode`\@=11
\def\matrice#1{\null \,\vcenter {\normalbaselines \m@th
\ialign {\hfil $##$\hfil &&\  \hfil $##$\hfil\crcr
\mathstrut \crcr \noalign {\kern -\baselineskip } #1\crcr
\mathstrut \crcr \noalign {\kern -\baselineskip }}}\,}

\def\petitematrice#1{\left(\null\vcenter {\normalbaselines \m@th
\ialign {\hfil $##$\hfil 
&&\thinspace  \hfil $##$\hfil\crcr
\mathstrut \crcr \noalign {\kern -\baselineskip } #1\crcr
\mathstrut \crcr \noalign {\kern -\baselineskip }}}\right)}

\catcode`\@=12

\def\qed{\quad\raise -2pt\hbox{\vrule\vbox to 10pt{\hrule width 4pt
   \vfill\hrule}\vrule}}

\def\cqfd{\penalty 500 \hbox{\qed}\par\smallskip}


%


\def\il{\bigl]\kern-.25em\bigl]}
\def\ir{\bigr]\kern-.25em\bigr]}

\def\iil{\bigl>\kern-.25em\bigl>}
\def\iir{\bigr>\kern-.25em\bigr>}



\def\Grille{\setbox13=\vbox to 5mm{\hrule width 110mm\vfill}
\setbox13=\vbox{\offinterlineskip
   \copy13\copy13\copy13\copy13\copy13\copy13\copy13\copy13
   \copy13\copy13\copy13\copy13\box13\hrule width 110mm}
\setbox14=\hbox to 5mm{\vrule height 65mm\hfill}
\setbox14=\hbox{\copy14\copy14\copy14\copy14\copy14\copy14
   \copy14\copy14\copy14\copy14\copy14\copy14\copy14\copy14
   \copy14\copy14\copy14\copy14\copy14\copy14\copy14\copy14\box14}
\ht14=0pt\dp14=0pt\wd14=0pt
\setbox13=\vbox to 0pt{\vss\box13\offinterlineskip\box14}
\wd13=0pt\box13}


\def\fleche(#1,#2)\dir(#3,#4)\long#5{%
\noalign{\nointerlineskip\leftput(#1,#2){\vector(#3,#4){#5}}\nointerlineskip}}


\def\hfl#1#2#3{\smash{\mathop{\hbox to#3{\rightarrowfill}}\limits
^{\scriptstyle#1}_{\scriptstyle#2}}}

\def\gfl#1#2#3{\smash{\mathop{\hbox to#3{\leftarrowfill}}\limits
^{\scriptstyle#1}_{\scriptstyle#2}}}


 \message{`lline' & `vector' macros from LaTeX}
 \catcode`@=11
\def\{{\relax\ifmmode\lbrace\else$\lbrace$\fi}
\def\}{\relax\ifmmode\rbrace\else$\rbrace$\fi}
\def\newcount{\alloc@0\count\countdef\insc@unt}
\def\newdimen{\alloc@1\dimen\dimendef\insc@unt}
\def\newwrite{\alloc@7\write\chardef\sixt@@n}

\newwrite\@unused
\newcount\@tempcnta
\newcount\@tempcntb
\newdimen\@tempdima
\newdimen\@tempdimb
\newbox\@tempboxa

\def\@spaces{\space\space\space\space}
\def\@whilenoop#1{}
\def\@whiledim#1\do #2{\ifdim #1\relax#2\@iwhiledim{#1\relax#2}\fi}
\def\@iwhiledim#1{\ifdim #1\let\@nextwhile=\@iwhiledim
        \else\let\@nextwhile=\@whilenoop\fi\@nextwhile{#1}}
\def\@badlinearg{\@latexerr{Bad \string\line\space or \string\vector
   \space argument}}
\def\@latexerr#1#2{\begingroup
\edef\@tempc{#2}\expandafter\errhelp\expandafter{\@tempc}%
\def\@eha{Your command was ignored.
^^JType \space I <command> <return> \space to replace it
  with another command,^^Jor \space <return> \space to continue without it.}
\def\@ehb{You've lost some text. \space \@ehc}
\def\@ehc{Try typing \space <return>
  \space to proceed.^^JIf that doesn't work, type \space X <return> \space to
  quit.}
\def\@ehd{You're in trouble here.  \space\@ehc}

\typeout{LaTeX error. \space See LaTeX manual for explanation.^^J
 \space\@spaces\@spaces\@spaces Type \space H <return> \space for
 immediate help.}\errmessage{#1}\endgroup}
\def\typeout#1{{\let\protect\string\immediate\write\@unused{#1}}}

\font\tenln    = line10
\font\tenlnw   = linew10

\newdimen\@wholewidth
\newdimen\@halfwidth
\newdimen\unitlength 

\unitlength =1pt


\def\thinlines{\let\@linefnt\tenln \let\@circlefnt\tencirc
  \@wholewidth\fontdimen8\tenln \@halfwidth .5\@wholewidth}
\def\thicklines{\let\@linefnt\tenlnw \let\@circlefnt\tencircw
  \@wholewidth\fontdimen8\tenlnw \@halfwidth .5\@wholewidth}

\def\linethickness#1{\@wholewidth #1\relax \@halfwidth .5\@wholewidth}

\newif\if@negarg

\def\lline(#1,#2)#3{\@xarg #1\relax \@yarg #2\relax
\@linelen=#3\unitlength
\ifnum\@xarg =0 \@vline
  \else \ifnum\@yarg =0 \@hline \else \@sline\fi
\fi}

\def\@sline{\ifnum\@xarg< 0 \@negargtrue \@xarg -\@xarg \@yyarg -\@yarg
  \else \@negargfalse \@yyarg \@yarg \fi
\ifnum \@yyarg >0 \@tempcnta\@yyarg \else \@tempcnta -\@yyarg \fi
\ifnum\@tempcnta>6 \@badlinearg\@tempcnta0 \fi
\setbox\@linechar\hbox{\@linefnt\@getlinechar(\@xarg,\@yyarg)}%
\ifnum \@yarg >0 \let\@upordown\raise \@clnht\z@
   \else\let\@upordown\lower \@clnht \ht\@linechar\fi
\@clnwd=\wd\@linechar
\if@negarg \hskip -\wd\@linechar \def\@tempa{\hskip -2\wd\@linechar}\else
     \let\@tempa\relax \fi
\@whiledim \@clnwd <\@linelen \do
  {\@upordown\@clnht\copy\@linechar
   \@tempa
   \advance\@clnht \ht\@linechar
   \advance\@clnwd \wd\@linechar}%
\advance\@clnht -\ht\@linechar
\advance\@clnwd -\wd\@linechar
\@tempdima\@linelen\advance\@tempdima -\@clnwd
\@tempdimb\@tempdima\advance\@tempdimb -\wd\@linechar
\if@negarg \hskip -\@tempdimb \else \hskip \@tempdimb \fi
\multiply\@tempdima \@m
\@tempcnta \@tempdima \@tempdima \wd\@linechar \divide\@tempcnta \@tempdima
\@tempdima \ht\@linechar \multiply\@tempdima \@tempcnta
\divide\@tempdima \@m
\advance\@clnht \@tempdima
\ifdim \@linelen <\wd\@linechar
   \hskip \wd\@linechar
  \else\@upordown\@clnht\copy\@linechar\fi}

\def\@hline{\ifnum \@xarg <0 \hskip -\@linelen \fi
\vrule height \@halfwidth depth \@halfwidth width \@linelen
\ifnum \@xarg <0 \hskip -\@linelen \fi}

\def\@getlinechar(#1,#2){\@tempcnta#1\relax\multiply\@tempcnta 8
\advance\@tempcnta -9 \ifnum #2>0 \advance\@tempcnta #2\relax\else
\advance\@tempcnta -#2\relax\advance\@tempcnta 64 \fi
\char\@tempcnta}

\def\vector(#1,#2)#3{\@xarg #1\relax \@yarg #2\relax
\@linelen=#3\unitlength
\ifnum\@xarg =0 \@vvector
  \else \ifnum\@yarg =0 \@hvector \else \@svector\fi
\fi}

\def\@hvector{\@hline\hbox to 0pt{\@linefnt
\ifnum \@xarg <0 \@getlarrow(1,0)\hss\else
    \hss\@getrarrow(1,0)\fi}}

\def\@vvector{\ifnum \@yarg <0 \@downvector \else \@upvector \fi}

\def\@svector{\@sline
\@tempcnta\@yarg \ifnum\@tempcnta <0 \@tempcnta=-\@tempcnta\fi
\ifnum\@tempcnta <5
  \hskip -\wd\@linechar
  \@upordown\@clnht \hbox{\@linefnt  \if@negarg
  \@getlarrow(\@xarg,\@yyarg) \else \@getrarrow(\@xarg,\@yyarg) \fi}%
\else\@badlinearg\fi}

\def\@getlarrow(#1,#2){\ifnum #2 =\z@ \@tempcnta='33\else
\@tempcnta=#1\relax\multiply\@tempcnta \sixt@@n \advance\@tempcnta
-9 \@tempcntb=#2\relax\multiply\@tempcntb \tw@
\ifnum \@tempcntb >0 \advance\@tempcnta \@tempcntb\relax
\else\advance\@tempcnta -\@tempcntb\advance\@tempcnta 64
\fi\fi\char\@tempcnta}

\def\@getrarrow(#1,#2){\@tempcntb=#2\relax
\ifnum\@tempcntb < 0 \@tempcntb=-\@tempcntb\relax\fi
\ifcase \@tempcntb\relax \@tempcnta='55 \or
\ifnum #1<3 \@tempcnta=#1\relax\multiply\@tempcnta
24 \advance\@tempcnta -6 \else \ifnum #1=3 \@tempcnta=49
\else\@tempcnta=58 \fi\fi\or
\ifnum #1<3 \@tempcnta=#1\relax\multiply\@tempcnta
24 \advance\@tempcnta -3 \else \@tempcnta=51\fi\or
\@tempcnta=#1\relax\multiply\@tempcnta
\sixt@@n \advance\@tempcnta -\tw@ \else
\@tempcnta=#1\relax\multiply\@tempcnta
\sixt@@n \advance\@tempcnta 7 \fi\ifnum #2<0 \advance\@tempcnta 64 \fi
\char\@tempcnta}

\def\@vline{\ifnum \@yarg <0 \@downline \else \@upline\fi}

\def\@upline{\hbox to \z@{\hskip -\@halfwidth \vrule
  width \@wholewidth height \@linelen depth \z@\hss}}

\def\@downline{\hbox to \z@{\hskip -\@halfwidth \vrule
  width \@wholewidth height \z@ depth \@linelen \hss}}

\def\@upvector{\@upline\setbox\@tempboxa\hbox{\@linefnt\char'66}\raise
     \@linelen \hbox to\z@{\lower \ht\@tempboxa\box\@tempboxa\hss}}

\def\@downvector{\@downline\lower \@linelen
      \hbox to \z@{\@linefnt\char'77\hss}}

\thinlines

\newcount\@xarg
\newcount\@yarg
\newcount\@yyarg
\newcount\@multicnt
\newdimen\@xdim
\newdimen\@ydim
\newbox\@linechar
\newdimen\@linelen
\newdimen\@clnwd
\newdimen\@clnht
\newdimen\@dashdim
\newbox\@dashbox
\newcount\@dashcnt
 \catcode`@=12


\newbox\tbox
\newbox\tboxa

\def\leftzer#1{\setbox\tbox=\hbox to 0pt{#1\hss}%
     \ht\tbox=0pt \dp\tbox=0pt \box\tbox}

\def\rightzer#1{\setbox\tbox=\hbox to 0pt{\hss #1}%
     \ht\tbox=0pt \dp\tbox=0pt \box\tbox}

\def\centerzer#1{\setbox\tbox=\hbox to 0pt{\hss #1\hss}%
     \ht\tbox=0pt \dp\tbox=0pt \box\tbox}

%
\def\image(#1,#2)#3{\vbox to #1{\offinterlineskip
    \vss #3 \vskip #2}}


\def\leftput(#1,#2)#3{\setbox\tboxa=\hbox{%
    \kern #1\unitlength
    \raise #2\unitlength\hbox{\leftzer{#3}}}%
    \ht\tboxa=0pt \wd\tboxa=0pt \dp\tboxa=0pt\box\tboxa}

\def\rightput(#1,#2)#3{\setbox\tboxa=\hbox{%
    \kern #1\unitlength
    \raise #2\unitlength\hbox{\rightzer{#3}}}%
    \ht\tboxa=0pt \wd\tboxa=0pt \dp\tboxa=0pt\box\tboxa}

\def\centerput(#1,#2)#3{\setbox\tboxa=\hbox{%
    \kern #1\unitlength
    \raise #2\unitlength\hbox{\centerzer{#3}}}%
    \ht\tboxa=0pt \wd\tboxa=0pt \dp\tboxa=0pt\box\tboxa}

\unitlength=1mm

\def\put(#1,#2)#3{\noalign{\nointerlineskip
                               \centerput(#1,#2){$#3$}
                                \nointerlineskip}}
\def\segment(#1,#2)\dir(#3,#4)\long#5{%
\leftput(#1,#2){\lline(#3,#4){#5}}}

\def\arbrea#1#2#3#4{\mathop{\hskip14pt
\vbox{\vskip1cm\offinterlineskip
\segment(0,0)\dir(0,1)\long{3}
\segment(0,3)\dir(1,1)\long{4}
\segment(0,3)\dir(-1,1)\long{4}
\rightput(-4,4.5){\hbox{$\scriptstyle#1$}}
\leftput(3.3,6.9){\hbox{$\scriptstyle#2$}}
\leftput(1,2){\hbox{$\scriptstyle#3$}}
\centerput(-4,7.5){\hbox{$\scriptstyle#4$}}
}\hskip14pt}\nolimits}

\def\carre{\mathop{\hbox{\kern5pt \vbox{
\offinterlineskip\segment(-1.5,-.5)\dir(0,1)\long{3}
\segment(-1.5,-0.5)\dir(1,0)\long{3}
\segment(1.5,-0.5)\dir(0,1)\long{3}
\segment(1.5,2.5)\dir(-1,0)\long{3}}\kern5pt }}}

\def\carresec{\mathop{\hbox{\kern5pt \vbox{
\offinterlineskip\segment(-1.5,-.5)\dir(0,1)\long{3}
\segment(1.8,-.8)\dir(-1,1)\long{3.6}
\segment(-1.5,-0.5)\dir(1,0)\long{3}
\segment(1.5,-0.5)\dir(0,1)\long{3}
\segment(1.5,2.5)\dir(-1,0)\long{3}}\kern5pt }}}

\def\rondsec{{{\bigcirc}\kern-7.5pt\backslash}\,}

\def\rond{\mathop{\raise.5pt\hbox{$\bigcirc$}}}

\def\rondbullet{\mathop{\hbox{\kern5pt \vbox{
\offinterlineskip
\centerput(0,2){\hbox{$\scriptstyle\bullet$}}
\centerput(0,0){\hbox{$\bigcirc$}}
\centerput(0,4){$\scriptstyle m$}
}}}}

\def\triang{\mathop{\raise.5pt\hbox{$\bigtriangledown$}}}
\def\triangup{\mathop{\raise.5pt\hbox{$\bigtriangleup$}}}

\def\triangbullet{\mathop{\hbox{\kern5pt \vbox{
\offinterlineskip
\centerput(0,2){\hbox{$\triang$}}
\centerput(1.7,3){$\scriptstyle\bullet$}
\centerput(3,4.4){$\scriptstyle m$}
}}}}

\def\triangsec{\mathop{\raise.5pt\hbox{$\textstyle\bigtriangledown\kern-6pt\backslash\,$}}}

\def\arbrebb#1#2#3#4#5#6#7#8{\mathop{\hskip14pt
\vbox{\vskip1.2cm\offinterlineskip
\segment(0,0)\dir(0,1)\long{3}
\segment(0,3)\dir(1,1)\long{4}
\segment(0,3)\dir(-1,1)\long{8}
\segment(-4.5,7)\dir(1,1)\long{4}
\leftput(-10,11.5){$\scriptstyle#1$}
\leftput(2.5,7.6){\hbox{$\scriptstyle#2$}}
\leftput(-1.3,10.8){\hbox{$\scriptstyle#3$}}
\rightput(-4,4.5){\hbox{$\scriptstyle#4$}}
\rightput(-8,9){\hbox{$\scriptstyle#5$}}
\leftput(0.2,11){\hbox{$\scriptstyle#6$}}
\leftput(0.5,1.5){\hbox{$\scriptstyle#7$}}
\leftput(4,5){\hbox{$\scriptstyle#8$}}
}\hskip14pt}\nolimits}

\def\arbrebbc#1#2#3#4#5#6#7#8#9{\mathop{\hskip14pt
\vbox{\vskip1.6cm\offinterlineskip
\segment(4,-1)\dir(0,1)\long{3}
\segment(4,2)\dir(1,1)\long{4}
\segment(0,6)\dir(1,1)\long{4}
\segment(0,6)\dir(-1,1)\long{8}
\segment(0,6)\dir(1,-1)\long{4}
\segment(-4.5,10)\dir(1,1)\long{4}
\leftput(-10,14.5){$\scriptstyle#1$}
\leftput(2.5,10.6){\hbox{$\scriptstyle#2$}}
\leftput(-1.3,13.8){\hbox{$\scriptstyle#3$}}
\rightput(-4,7.5){\hbox{$\scriptstyle#4$}}
\rightput(-8,12){\hbox{$\scriptstyle#5$}}
\leftput(0.2,14){\hbox{$\scriptstyle#6$}}
\rightput(-0.5,4.5){\hbox{$\scriptstyle#7$}}
\leftput(4,8){\hbox{$\scriptstyle#8$}}
\centerput(8,6.5){$\scriptstyle#9$}
}\hskip14pt}\nolimits}

\def\arbrebbz#1#2#3#4#5#6#7#8#9{\mathop{\hskip14pt
\vbox{\vskip2cm\offinterlineskip
\segment(0,6)\dir(1,1)\long{10}
\segment(6,12)\dir(-1,1)\long{4}
\segment(0,6)\dir(-1,1)\long{8}
\segment(0,6)\dir(0,-1)\long{4}
\segment(-4.5,10)\dir(1,1)\long{4}
\leftput(-10,14.5){$\scriptstyle#1$}
\centerput(6,9){\hbox{$\scriptstyle#2$}}
\leftput(-1.3,13.5){\hbox{$\scriptstyle#3$}}
\rightput(-4,7.5){\hbox{$\scriptstyle#4$}}
\rightput(1,15){\hbox{$\scriptstyle#5$}}
\leftput(1,16.5){\hbox{$\scriptstyle#6$}}
\rightput(2,3.5){\hbox{$\scriptstyle#7$}}
\leftput(3.5,9.5){\hbox{$\scriptstyle#8$}}
\centerput(10,16.5){$\scriptstyle#9$}
}\hskip14pt}\nolimits}

\def\arbrebb#1#2#3#4#5#6#7#8{\mathop{\hskip14pt
\vbox{\vskip1.2cm\offinterlineskip
\segment(0,0)\dir(0,1)\long{3}
\segment(0,3)\dir(1,1)\long{4}
\segment(0,3)\dir(-1,1)\long{8}
\segment(-4.5,7)\dir(1,1)\long{4}
\leftput(-10,11.5){$\scriptstyle#1$}
\leftput(2.5,7.6){\hbox{$\scriptstyle#2$}}
\leftput(-1.3,10.8){\hbox{$\scriptstyle#3$}}
\rightput(-4,4.5){\hbox{$\scriptstyle#4$}}
\rightput(-8,9){\hbox{$\scriptstyle#5$}}
\leftput(0.2,11){\hbox{$\scriptstyle#6$}}
\leftput(0.5,1.5){\hbox{$\scriptstyle#7$}}
\leftput(4,5){\hbox{$\scriptstyle#8$}}
}\hskip14pt}\nolimits}

\def\arbrebbj#1#2#3#4#5#6#7#8{\mathop{\hskip14pt
\vbox{\vskip1.2cm\offinterlineskip
\segment(0,0)\dir(0,1)\long{3}
\segment(0,3)\dir(-1,1)\long{4}
\segment(0,3)\dir(1,1)\long{8}
\segment(4.5,7)\dir(-1,1)\long{4}
\leftput(-10,11.5){$\scriptstyle#1$}
\leftput(7.5,9){\hbox{$\scriptstyle#2$}}
\leftput(-1.3,10.8){\hbox{$\scriptstyle#3$}}
\centerput(-4,6.5){\hbox{$\scriptstyle#4$}}
\rightput(10,11.5){\hbox{$\scriptstyle#5$}}
\leftput(-1.5,11.5){\hbox{$\scriptstyle#6$}}
\rightput(-4,4.5){\hbox{$\scriptstyle#7$}}
\leftput(4,5){\hbox{$\scriptstyle#8$}}
}\hskip14pt}\nolimits}

\def\Deltaa{\mathop{\hbox{$\Delta$}}\limits}
\def\brullet{{\scriptscriptstyle\bullet}}

\def\eoc{\mathop{\rm eoc}\limits}

\def\pom{\mathop{\rm pom}\limits}

\auteurcourant={DOMINIQUE FOATA AND GUO-NIU HAN}
\titrecourant={TREE SECANT CALCULUS}


\rightline{March 8, 2013}
\bigskip
\centerline{\bf 
Secant Tree Calculus}
\bigskip
\centerline{\sl Dominique Foata and Guo-Niu Han}
\footnote{}{
}

\bigskip
{\narrower\narrower
\eightpoint
\noindent
{\bf Abstract}.\quad
A true Tree Calculus is being developed to make a joint study of the two statistics ``eoc'' (end of minimal chain) and ``pom'' (parent of maximum leaf) on the set of secant trees. Their joint distribution restricted to the set $\{\eoc-\pom\le 1\}$ is shown to satisfy two partial difference equation systems, to be symmetric and to be expressed in the form of an explicit three-variable generating function.

}
\bigskip
\centerline{\bf 1. Introduction}

\medskip
As was done in our previous paper, whose purpose was to 
evaluate the distribution of 
a two-variable statistic defined on the set of {\it tangent trees} [FH13],
we pursue the same goal with the secant trees, by using an appropriate
{\it Tree Calculus}.
The secant trees, defined in full detail in the sequel, are labeled binary trees; they are called {\it secant}, because the number of them with $2n$ nodes is equal to the {\it secant} number $E_{2n}$, namely, the coefficient of $u^{2n}/(2n)!$ in the Taylor expansion of $\sec u$:
$$
\leqalignno{\noalign{\vskip-6pt}
\qquad\sec u={1\over \cos u}&=\sum_{n\ge 0}
{u^{2n}\over (2n)!}E_{2n}&(1.1)\cr
\noalign{\vskip-5pt}
&=1+{u^2\over 2!}1+{u^4\over 4!}5+{u^6\over 6!}61+{u^8\over
8!}1385+ {u^{10}\over 10!}50521+\cdots\cr 
}
$$
(see, e.g., [Ni23, p. 177-178], [Co74, p. 258-259]).
The purpose of this paper is to calculate the {\it joint distribution} of the two statistics `eoc'' (end of minimal chain) and ``pom'' (parent of maximum leaf) on the set of secant trees. To achieve this and in particular derive an explicit  three-variable generating function two {\it partial difference equation systems} are to be solved. The solutions of those systems are based on {\it Tree Calculus}, which consists of partitioning each set of secant trees into
smaller subsets, and developing a natural algebra on those subsets for solving the difference equations.

\smallskip
1.1. {\it The secant trees}.\quad 
The trees in question form a subclass of the so-called {\it increasing, binary} trees. The latter trees having $n$ nodes labeled $1,2,\ldots,n$ are in one-to-one correspondence with the $n!$ ordinary permutations of the sequence $12\cdots n$. See, for instance,  Viennot [Vi88, chap.~3]. Their actual definition is next stated by using the traditional vocabulary on trees, such as node, leaf, child, root, \dots\  In particular, when a node is not a leaf, it is said to be an {\it interior node}.

\medskip\smallskip
{\it Definition}.\quad
For each positive integer~$n$ an {\it increasing} tree of size~$n$ is defined by the following axioms:

(1) it is a {\it labeled\/} tree with $n$ nodes, labeled
$1,2,\ldots, n$; the node labeled~1 being the {\it root};

(2) each node has no child (then called a {\it leaf\/}), or
one child, or two children;

(3) the label of each node is smaller than the label of its children, if any;

(4) the tree is planar and each child of each node is, either on the
left (it is then called the {\it left child\/}), or on the right (the
{\it right child\/}); moreover, the tree can be embedded on the Euclidean plane as follows: the root has coordinates $(0,0)$, the left child (if any) $(-1,1)$, the right child (if any) $(1,1)$, the grandchildren (if any) $(-3/2,2)$, $(-1/2,2)$, $(1/2,2)$, $(3/2,2)$, the greatgrandchildren (if any) $(-7/4,3)$, $(-5/4,3)$, \dots~, $(7/4,3)$, etc. With this convention all the nodes have different abscissas. The node having the maximum abscissa is then defined in a unique manner. Call it the {\it rightmost node}; it is either a leaf, or a node having a left child, but no right child.

Consider the orthogonal projections of those $n$ nodes onto a horizontal axis. Reading the labels of those projected $n$ nodes from left to right yields a permutation $\sigma=\sigma(1)\sigma(2)\cdots \sigma(n)$ of $12\cdots n$. This {\it projection} defines a bijection of the set of all increasing trees with~$n$ nodes onto the permutation group~${\goth S}_{n}$.

The trees which correspond to permutations $x_{1}x_{2}\cdots x_{n}$ having the property that $x_{1}>x_{2}$, $x_{2}<x_{3}$, $x_{3}>x_{4}$, \dots~, in an alternating way, usually named {\it alternating}, are called {\it complete increasing}. In an equivalent manner, an increasing tree with~$n$ node is said to be {\it complete} (see [Vi88],  [KPP94], [St99]), if axioms (1)--(4) hold with the further property

(5) every node is, either a leaf, or a node with two children, except 
the rightmost node, which has one {\it left} child, but no right child when~$n$ is {\it even}. This rightmost node will then be referred to as being the {\it one-child node\/}.

\medskip\smallskip
Each complete increasing tree with~$n$ nodes is simply called {\it secant} (resp. {\it tangent\/}) whenever $n$ is {\it even} (resp. {\it odd\/}). For each $n\ge 1$ let~${\goth T}_{n}$ be the set of all complete increasing trees of size~$n$. In Fig.~1.1 the five complete increasing trees from~${\goth T}_{4}$, accordingly the five {\it secant} trees, have been drawn, together with the projections of their node labels on the horizontal axis.  
Notice that each of those projections, when read from left to right, forms an {\it alternating permutation} of the sequence $1\,2\,3\,4$. Under each tree have been calculated the two statistics ``eoc'' and ``pom'' defined below.

\goodbreak

\newbox\boxarbre
\newbox\boxarbrea
\newbox\boxarbreb
\newbox\boxarbrec
\newbox\boxarbred

\setbox\boxarbre=\vbox{\vskip
17mm\offinterlineskip 
\centerput(10,16){$3$}
\centerput(0,11){$4$}
\centerput(5,2){$1$}
\centerput(15,7){$2$}
\segment(0,10)\dir(1,-1)\long{5}
\segment(5,5)\dir(2,1)\long{10}
\segment(15,10)\dir(-1,1)\long{5}
\segment(0,0)\dir(1,0)\long{15}
\centerput(0,-4){$4$}
\centerput(5,-4){$1$}
\centerput(10,-4){$3$}
\centerput(15,-4){$2$}
\centerput(2,-9){$\eoc={}$}
\centerput(2,-13){$\pom={}$}
\centerput(12.5,-9){$3$}
\centerput(12.5,-13){$1$}
}

\setbox\boxarbrea=\vbox{\vskip
17mm\offinterlineskip 
\centerput(10,16){$4$}
\centerput(0,11){$3$}
\centerput(5,2){$1$}
\centerput(15,7){$2$}
\segment(0,10)\dir(1,-1)\long{5}
\segment(5,5)\dir(2,1)\long{10}
\segment(15,10)\dir(-1,1)\long{5}
\segment(0,0)\dir(1,0)\long{15}
\centerput(0,-4){$3$}
\centerput(5,-4){$1$}
\centerput(10,-4){$4$}
\centerput(15,-4){$2$}
\centerput(7.5,-9){$4$}
\centerput(7.5,-13){$2$}
}

\setbox\boxarbreb=\vbox{\vskip
17mm\offinterlineskip 
\centerput(10,16){$4$}
\centerput(0,11){$2$}
\centerput(5,2){$1$}
\centerput(15,7){$3$}
\segment(0,10)\dir(1,-1)\long{5}
\segment(5,5)\dir(2,1)\long{10}
\segment(15,10)\dir(-1,1)\long{5}
\segment(0,0)\dir(1,0)\long{15}
\centerput(0,-4){$2$}
\centerput(5,-4){$1$}
\centerput(10,-4){$4$}
\centerput(15,-4){$3$}
\centerput(7.5,-9){$3$}
\centerput(7.5,-13){$2$}
}

\setbox\boxarbreb=\vbox{\vskip
17mm\offinterlineskip 
\centerput(10,16){$4$}
\centerput(0,11){$2$}
\centerput(5,2){$1$}
\centerput(15,7){$3$}
\segment(0,10)\dir(1,-1)\long{5}
\segment(5,5)\dir(2,1)\long{10}
\segment(15,10)\dir(-1,1)\long{5}
\segment(0,0)\dir(1,0)\long{15}
\centerput(0,-4){$2$}
\centerput(5,-4){$1$}
\centerput(10,-4){$4$}
\centerput(15,-4){$3$}
\centerput(7.5,-9){$3$}
\centerput(7.5,-13){$2$}
}

\setbox\boxarbrec=\vbox{\vskip
17mm\offinterlineskip 
\centerput(5,8){$2$}
\centerput(0,11.5){$4$}
\centerput(15,6){$1$}
\centerput(10,11){$3$}
\segment(0,15)\dir(3,-2)\long{15}
\segment(5,11.5)\dir(2,1)\long{5}
\segment(0,0)\dir(1,0)\long{15}
\centerput(0,-4){$4$}
\centerput(5,-4){$2$}
\centerput(10,-4){$3$}
\centerput(15,-4){$1$}
\centerput(7.5,-9){$3$}
\centerput(7.5,-13){$2$}
}

\setbox\boxarbred=\vbox{\vskip
17mm\offinterlineskip 
\centerput(5,8){$2$}
\centerput(0,11.5){$3$}
\centerput(15,6){$1$}
\centerput(10,11){$4$}
\segment(0,15)\dir(3,-2)\long{15}
\segment(5,11.5)\dir(2,1)\long{5}
\segment(0,0)\dir(1,0)\long{15}
\centerput(0,-4){$3$}
\centerput(5,-4){$2$}
\centerput(10,-4){$4$}
\centerput(15,-4){$1$}
\centerput(7.5,-9){$2$}
\centerput(7.5,-13){$3$}
}

\vskip-12pt

$$
\hskip4mm
\box\boxarbre\hskip22mm\box\boxarbrea
\hskip22mm\box\boxarbrec\hskip22mm\box\boxarbred
\hskip22mm\box\boxarbreb\hskip2cm
$$

\bigskip\bigskip\bigskip
\centerline{Fig. 1.1}

\medskip
1.2. {\it The main two statistics}.\quad
Let $t\in {\goth T}_{n}$ $(n\ge 1)$. 
If a node labeled~$a$ has two children labeled $b$ and~$c$, define $\min a:=\min\{b,c\}$; if it has one child~$b$, let $\min a:=b$. The {\it minimal chain} of~$t$ is defined to be the sequence
$a_1\rightarrow a_2\rightarrow a_3\rightarrow\cdots
\rightarrow a_{j-1}\rightarrow a_j$, with the following properties: (i) $a_1=1$ is the label of the root; (ii) for
$i=1,2,\ldots,j-2$ the $(i+1)$-st term~$a_{i+1}$ is the label of an interior node and $a_{i+1}=\min a_i$; (iii) $a_j$ is the node of a leaf. Define the ``{\bf e}nd {\bf o}f the minimal {\bf c}hain'' of~$t$ to be $\eoc(t):=a_{j}$. If the leaf with the maximum label~$n$ is incident to a node labeled~$k$, define its ``{\bf p}arent {\bf o}f the {\bf m}aximum leaf'' to be $\pom(t):=k$. See Fig.~1.1. 

Those two statistics have been introduced by Poupard [Po88] on {\it tangent} trees ($n$ even). She proved that ``eoc'' and ``1+pom''
were equidistributed on each set ${\goth T}_{2n+1}$ of tangent trees, also that their common univariable distribution satisfied a finite difference equation system; she further calculated their generating function. In [FH13] it was proved that the equidistribution actually holds on every set ${\goth T}_{n}$ $(n\ge 1)$ by constructing an explict bijection~$\phi$ of ${\goth T}_{n}$ onto itself with the property that:
$1\!+\!\pom(t)=\eoc\phi(t)$ for all~$t$.

Working with alternating permutations it was shown in [FH12]
that the distribution of each statistic ``eoc,'' ``pom'' on {\it secant\/} trees satisfied the {\it same} finite difference system, as introduced by Poupard for tangent trees, but the initial conditions were different. Their generating function could also be calculated.

The next step was to study the {\it joint} distribution of the pair 
$(\eoc,\pom)$ on each set ${\goth T}_{n}$, that is to say, letting
$$f_{n}(m,k):=\#\{t\in {\goth T}_{n}\!:\! \eoc (t)=m\ {\rm and}\ \pom (t)=k\}\leqno(1.2)
$$
see whether those numbers are solutions of a {\it partial} finite difference equation system, and try to calculate their generating function. This program has been achieved in our paper [FH13] for the sets of {\it tangent} trees ${\goth T}_{2n+1}$. The purpose of this paper is then to pursue this program for {\it secant} trees.

\goodbreak
1.3. {\it The joint distribution}.\quad
The first values of the matrices $M_{2n}:= (f_{2n}(m,k))$ 
$(2\le m\le 2n; 1\le k\le 2n-1)$ are listed in Table 1.1. 
For each matrix have been evaluated the {\it row sums} $f_{2n}(m,\brullet):=\sum_{k}f_{2n}(m,k)
=\#\{t\in {\goth T}_{2n}\!:\! \eoc (t)=m\}$
(resp. {\it column sums} $f_{2n}(\brullet,k):=\sum_{m}f_{2n}(m,k)=\#\{t\in {\goth T}_{2n}\!:\! \pom (t)=k\}$) on the rightmost column (resp. the bottom row). In the South-East-corner is written the total sum $f_{2n}(\brullet,\brullet):=\sum_{m,k}f_{2n}(m,k)$, equal to the {\it secant number} $E_{2n}$. Keeping in mind that alternating permutations are equidistributed with complete increasing trees, the identity
$$f_{2n}(\brullet,\brullet)=\sum_{m,k}f_{2n}(m,k)=E_{2n}\leqno(1.3)$$ is a consequence of the old result derived by D\'esir\'e Andr\'e [An1979, An1881]. Also, as ``1+pom'' and ``eoc'' are equidistributed, we have:
$$
f_{2n}(\brullet,k-1)=f_{2n}(k,\brullet)\quad(2\le k\le 2n).
\leqno(1.4)
$$

\midinsert

\medskip
{\eightpoint
$$M_{2}=\vcenter{\offinterlineskip
\halign{\vrule height10pt width.4pt
depth4pt$\hfil\ #\ $&\vrule height10pt width.4pt
depth4pt$\hfil\ #\ $&%
&\vrule height10pt width.4pt depth4pt$\hfil\,#\,\hfil$%
\vrule height10pt width.4pt depth4pt\cr 
\noalign{\hrule}
k=&1&f_2(m,.)\cr 
\noalign{\hrule}
m=2&1&1\cr
\noalign{\hrule}
f_2(.,k)&1&E_2=1\cr
\noalign{\hrule}
}}\qquad M_{4}=\vcenter{\offinterlineskip
\halign{\vrule height10pt width.4pt
depth4pt$\hfil\ #\ $&\vrule height10pt width.4pt
depth4pt$\hfil\ #\ $&
$\hfil\ #\ $&$\hfil\ #\ $%
&\vrule height10pt width.4pt depth4pt$\hfil\,#\,\hfil$%
\vrule height10pt width.4pt depth4pt\cr 
\noalign{\hrule}
k=&1&2&3&f_4(m,.)\cr 
\noalign{\hrule}
 m=2& .  &   .  & 1 &1 \cr
 3& 1 &  2 &.&3\cr
4&  .  & 1 &.&1\cr
\noalign{\hrule}
f_4(.,k)&1&3&1&E_4=5\cr
\noalign{\hrule}
}}
  $$
  
$$M_{6}=\vcenter{\offinterlineskip
\halign{\vrule height10pt width.4pt
depth4pt$\hfil\ #\ $&\vrule height10pt width.4pt
depth4pt$\hfil\ #\ $&
$\hfil\ #\ $&$\hfil\ #\ $&$\hfil\ #\ $&$\hfil\ #\ $%
&\vrule height10pt width.4pt depth4pt$\hfil\,#\,\hfil$%
\vrule height10pt width.4pt depth4pt\cr 
\noalign{\hrule}
k=&1&2&3&4&5&f_6(m,.)\cr 
\noalign{\hrule}
 m=2   &    .  &   . &    1  &   3  &   1 &5    \cr
  3&      1  &   2  &   .  &   9  &   3 &15    \cr
   4&    3   &  7  &  10  &   .  &   1 &21    \cr
  5&  1  &   4  &   8  &   2  &   .  &15   \cr
6& . &        2  &   2  &   1  &   .  &5   \cr
\noalign{\hrule}
f_6(.,k)&5&15&21&15&5&E_6=61\cr
\noalign{\hrule}
}}$$

$$
M_{8}=\vcenter{\offinterlineskip
\halign{\vrule height10pt width.4pt
depth4pt$\hfil\ #\ $&\vrule height10pt width.4pt
depth4pt$\hfil\ #\ $&
$\hfil\ #\ $&$\hfil\ #\ $&$\hfil\ #\ $&$\hfil\ #\ $%
&$\hfil\ #\ $&$\hfil\ #\ $
&\vrule height10pt width.4pt depth4pt$\hfil\,#\,\hfil$%
\vrule height10pt width.4pt depth4pt\cr 
\noalign{\hrule}
k=&1&2&3&4&5&6&7&f_8(m,.)\cr 
\noalign{\hrule}
m=2 &       .  &   .  &   5  &  15 &   21 &   15 &    5 &61  \cr
 3&     5  &  10  &   . &   45 &   63 &   45 &   15&183  \cr
 4&    15  &  35  &  50  &   . &  101 &   63  &  21 &285\cr
 5&     21  &  54  &  86  & 106  &   . &   45  &  15 &327   \cr
  6&      15 &   46 &   82 &   87 &   50 &    .  &   5 &285   \cr
7&    5  &  22  &  46  &  60 &   40 &   10  &   .  &183  \cr
8&      . &   16  &  16 &   14  &  10   &  5  &   . &61   \cr
\noalign{\hrule}
f_8(.,k)&61&183&285&327&285&183&61&E_8=1385\cr
\noalign{\hrule}
}}$$

$$
M_{10}=\vcenter{\offinterlineskip
\halign{\vrule height10pt width.4pt
depth4pt$\hfil\, #\, $&\vrule height10pt width.4pt
depth4pt$\hfil\, #\, $&
$\hfil\, #\, $&$\hfil\, #\, $&$\hfil\, #\, $&$\hfil\, #\, $%
&$\hfil\, #\, $&$\hfil\, #\, $&$\hfil\, #\, $&$\hfil\, #\, $
&\vrule height10pt width.4pt depth4pt$\hfil\,#\,\hfil$%
\vrule height10pt width.4pt depth4pt\cr 
\noalign{\hrule}
k=&1&2&3&4&5&6&7&8&9&f_{10}(m,.)\cr 
\noalign{\hrule}
m=2&        .  &    . &   61  & 183 &  285  & 327 &  285 &  183 &   61&1385  \cr
 3&       61 &  122 &    . &  549 &  855 &  981 &  855 & 549  & 183& 4155  \cr      4&183& 427& 610& .&1405&1575&1341& 855& 285&6681\cr
 5& 285& 720&1132&1466& . & 1989&1575& 981& 327 &8475\cr
6& 327& 884&1460&1863&2050& .  & 1405& 855& 285& 9129\cr
7&285 &  836 & 1448 & 1838 & 1870 & 1466 &. & 549  & 183&8475  \cr
  8& 183 &  606 & 1110 & 1466 & 1490 & 1155 &  610  &  . &    61&6681  \cr
  9& 61 &  288  & 588  & 854  & 950 &  804 &  488  & 122 &.&4155\cr
 10&.  &   272  & 272 &  256  & 224 &  178 &  122  &  61&.&1385 \cr
 \noalign{\hrule}
 f_{10}(.,k)&1385&4155&6681&8475&9129&8475&6681&4155&1385&E_{10}=50521\cr
 \noalign{\hrule} 
}}
$$

}
\centerline{Table 1.1: the matrices $M_{2n}$ $(1\le n\le 5)$}

\endinsert

On each entry $f_{2n}(m,k)$ may be defined two {\it partial differences} with respect to~$m$ and~$k$ as follows:
$$
\leqalignno{\Deltaa_{m}f_{2n}(m,k)&:=f_{2n}(m+1,k)-f_{2n}(m,k);&(1.5)\cr
\Deltaa_{k}f_{2n}(m,k)&:=f_{2n}(m,k+1)-f_{2n}(m,k).&(1.6)\cr}
$$
By convention $f_{2n}(m,k):=0$ if $(m,k)\not\in [2,2n]\times[1,2n-1]$. Our main results are the following.

\proclaim Theorem 1.1.  The finite difference equation systems hold:
$$\leqalignno{\qquad\quad
{\Deltaa_{m}}^2f_{2n}(m,k)+4\,f_{2n-2}(m,k-2)&=0\ 
(2\le m\le k-3<k\le 2n-1);&(R\,1)\cr
{\Deltaa_{k}}^2f_{2n}(m,k)+4\,f_{2n-2}(m,k)&=0\ 
(2\le m\le k-1<k\le 2n-3).&(R\,2)\cr
}
$$

The first two top rows and rightmost two columns of the upper triangles $\{f_{2n}(m,k):2\le m<k\le 2n-1\}$ $(n\ge 2)$ can be evaluated in function of the row or column sums
$f_{2n-2}(\brullet,k)$, $f_{2n-2}(m,\brullet)$, as is now stated.

\proclaim Theorem 1.2. We have:
$$
\eqalignno{
	(1.7)\qquad f_{2n}(2,k)&=f_{2n-2}(\brullet,k-2)
	=f_{2n-2}(k-1,\brullet) \quad (3\le k\le 2n-1); \cr
&&\hbox{[First top row]}\cr
f_{2n}(3,k)&=3\,f_{2n}(2,k)\quad (4\le k\le 2n-1);&\hbox{[Second top row]}\cr 
f_{2n}(m,2n-1)&=f_{2n-2}(m,\brullet)=f_{2n-2}(\brullet,m-1)\quad (2\le m\le 2n-2);\cr
&&\hbox{[Rightmost column]}\cr
f_{2n}(m,2n-2)&=3\,f_{2n}(m,2n-1)
\quad(2\le m\le 2n-3).\cr
&&\hbox{[Next to rightmost column]}\cr
}$$

\proclaim Proposition 1.3. We further have:
$$
\leqalignno{
f_{2}(\brullet,1)=1;\quad
f_{2n}(\brullet,1)&=f_{2n-2}(\brullet,\brullet)=E_{2n-2}
\ (n\ge 2);&(1.8)\cr
f_{2n}(\brullet,2)&=3f_{2n-2}(\brullet,\brullet)=3E_{2n-2}\ (n\ge 2);&(1.9)\cr
{\Deltaa_{m}}^2\,f_{2n}(m,\brullet)+4\,f_{2n-2}(m,\brullet)&=0
\quad(2\le m\le 2n-2);&(R\,3)\cr
{\Deltaa_{k}}^2\,f_{2n}(\brullet,k)+4\,f_{2n-2}(\brullet,k)&=0
\quad(1\le k\le 2n-3).&(R\,4)\cr}
$$

\proclaim Theorem 1.4. The previous two theorems and Proposition~$1.3$ provide an explicit algorithm for calculating the entries of the {\it upper triangles} of the matrices $M_{2n}=(f_{2n}(m,k))$ $(n\ge 1)$.  

The entries of the {\it lower triangles} $\{f_{2n}(m,k):1\le k<m\le 2n\}$ in Tables~1.1 have been calculated directly by means of formula (1.2). Contrary to the {\it upper triangles} we do not have any explicit numerical algorithm to get them; only the entries situated on the three sides of those lower triangles can be directly evaluated, as shown in Section~6. See, in particular, Table~6.1. For the {\it upper triangles}, we can derive an explicit generating function, as stated in the next theorem.

\proclaim Theorem 1.5. The triple exponential generating function for the upper triangles of the matrices~$(f_{2n}(m,k))$ is given by
$$
\displaylines{(1.10)\quad
\sum_{2\le m<k\le 2n-1}
f_{2n}(m,k){x^{2n-k-1}\over (2n-k-1)!}
{y^{k-m-1}\over (k-m-1)!}
{z^{m-2}\over (m-2)!}\hfill\cr
\hfill{}=
{\cos(2y)+2\,\cos(2(x-z))-\cos(2(z+x))
\over 2\,\cos^3(x+y+z)}
.\quad\cr
}$$

The right-hand side of (1.10) is symmetric in $x, z$. Hence, the change $x\leftrightarrow z$ in the left-hand side of (1.10) shows that
$$
f_{2n}(2n+1-k,2n+1-m)=f_{2n}(m,k).\leqno(1.11)
$$
The upper triangles already mentioned are then symmetric with respect to their counter-diagonals. This result can also be extended as follows.

\proclaim Theorem 1.6. Let ${\rm Up}(2n)$ be the set of all $(m,k)$ 
from $[2,2n]\times[1,2n-1]$ such that either $m-1\le k$, $(m,k)=(3,1)$, or $(m,k)=(2n,2n-2)$. Then $(1.11)$ holds for every $(m,k)\in 
{\rm Up}(2n)$.

Those theorems and proposition are proved in the next sections, once the main ingredients on Tree Calculus have been developed, as done in the next section.

\vfill\eject
\centerline {\bf 2. Tree Calculus}

\medskip
We adopt the following notations and conventions: for each triple $(n,m,k)$ let ${\goth T}_{2n,m,k}$  (resp. ${\goth T}_{2n,m,\brullet}$, resp. ${\goth T}_{2n,\brullet,k}$) denote the subset of ${\goth T}_{2n}$ of all trees~$t$ such that $\eoc(t)=m$ and $\pom(t)=k$ (resp. $\eoc(t)=m$, resp. $\pom(t)=k$). 
Also, {\it symbols representing families of trees will also designate their cardinalities}. With this convention ${\goth T}_{n,m,k}:=\#{\goth T}_{n,m,k}$. The matrix of the upper triangles ${\goth T}_{2n,m,k}$ $(2\le m<k\le 2n-1)$ will be denoted by ${\rm Upper}({\goth T}_{2n})$. 
\medskip

Subtrees (possibly empty) are indicated by the symbols ``$\rond$," ``$\triang$", or ``$\carre$." The notation ``$\carresec\,$'' (resp. ``$\rondsec\,$,'' resp. ``$\triang\kern-8pt\backslash\,$'') is used to indicate that the subtree ``$\carre$'' (resp. ``$\rond$,'' resp. ``$\triang$'') {\it contains the one-child node or is empty}.
Letters occurring below or next to subtrees are labels of their roots. 
The end of the minimal chain in each tree is represented by a bullet ``$\bullet$.'' 

In the sequel certain families of secant trees will be represented by symbols, called {\it trunks}. For example, the symbol 
$$
A=\arbrea{}{\bullet\, m}b{\carre}
$$
is the trunk that designates the {\it family} of all trees~$t$ from the underlying set ${\goth T}_{2n}$ having a node $b$, parent of both a subtree  ``$\carre$" and the leaf~$m$, which is also the end of the minimal chain. 
Notice that unlike the secant trees which are ordered, the trunk is {\it unordered}, so that
$$
\arbrea{}{\bullet\, m}b{\carre}
\quad \equiv \quad
\arbrea{}{\carre}b{\lower2pt\hbox{${\scriptstyle\quad \bullet\;  m}$}}.
\leqno{(2.1)}
$$
When the subtree ``$\carre$" contains the one-child node or is empty, we let
$$
A(\carresec):=
\arbrea{}{\bullet\, m}b{\carresec}
\quad \equiv \quad
\arbrea{}{\carresec}b{\lower2pt\hbox{${\scriptstyle\quad \bullet\;  m}$}},
$$
be the family of all trees~$t$ from the underlying set ${\goth T}_{2n}$ having a node $b$,  parent of the {\it right} child  ``$\carresec$" and the {\it left} child leaf~$m$. 
Let $A(\star)$ be the set of all trees from $A$ such that 
``$\carre$" is not empty and does not contain the one-child node.
We have the following decomposition
$$
A=A(\star)+A(\carresec).\leqno(2.2)
$$

When a further condition $(C)$ is imposed on a trunk $A$ we shall use the 
notation $[A,(C)]$. 
For example, the symbol
$$ 
B=[\quad\arbrea{}{\bullet\, m}b{\carre}\quad,a\ ]
$$
is the trunk that
has the same characteristic of the trunk $A$ as above with the further property that the node labeled~$a$ belongs, {\it neither} to the subtree of root~$b$, {\it nor} to the path going from root~1 to~$b$. 

\medskip

In our Tree Calculus we shall mostly compare the cardinalities of certain pairs of trunks, as shown in the following two examples.

\medskip
{\it Example $1$}. The two trunks
$$
[\quad\arbrea{m+1}{\lower4.5pt\hbox{\kern-3pt$\rond\atop
m+2$}}{} {\lower2pt\hbox{$\bullet$}}\quad,m]\quad
\hbox{\ and \ }\quad[\quad\arbrea{m}{\kern-5pt \lower3pt\hbox{$\rond\atop \ m+1$}}{}{\lower2pt
\hbox{$\bullet$}}\quad,m+2]
$$
have the same cardinalities, by using the bijection 
${\ \ m\ \ m+1 \ m+2 \choose m+2 \ \ m \ \ m+1}$.

\medskip
{\it Example $2$}. To compare the cardinalities of the following two trunks
$$
C=\quad\arbrebbc{\rond}{\carre}{}{k+2}{}{2n}{k+1}{}{\triang}
\raise3pt\hbox{$\scriptstyle k$}
\hbox{\quad and \quad}
D = \quad\arbrebbc{\rond}{\lower2pt\hbox{$\;$}}{\raise2pt
\hbox{$\kern-2pt
\carre$}}{k+2}{}{}{k+1}{\raise5pt\hbox{$\scriptstyle\kern4pt
2n$}}{\triang}\raise3pt\hbox{$\scriptstyle k$}
{\raise10pt\hbox{$\scriptstyle\kern-8pt $}}\qquad,
$$
we decompose them according to the location of the
one-child node:
$$\leqalignno{
C
&=\quad\arbrebbc{\rond}{\carre}{}{k+2}{}{2n}{k+1}{}{\triang}
\raise3pt\hbox{$\scriptstyle k$} \quad
+
\quad \arbrebbc{\rond}{\carre}{}{k+2}{}{2n}{k+1}{}{\textstyle\triang\kern-8pt\backslash}
\raise3pt\hbox{$\scriptstyle k$} \quad
+
\quad \arbrebbc{\rond}{\carresec}{}{k+2}{}{2n}{k+1}{}{\triang}
\raise3pt\hbox{$\scriptstyle k$} \quad
+ 
\quad \arbrebbc{\textstyle\rondsec}{\carre}{}{k+2}{}{2n}{k+1}{}{\triang}
\raise3pt\hbox{$\scriptstyle k$}
\cr
&:=C(\star)+C(\triang\kern-8pt\backslash\,)+C(\carresec)+C(\rondsec);\cr
\noalign{\bigskip}
D&=
\quad \arbrebbc{\rond}{\lower2pt\hbox{$\;$}}{\raise2pt
\hbox{$\kern-2pt
\carre$}}{k+2}{}{}{k+1}{\raise5pt\hbox{$\scriptstyle\kern4pt
2n$}}{\triang}\raise3pt\hbox{$\scriptstyle k$}
{\raise10pt\hbox{$\scriptstyle\kern-8pt $}}\qquad
+\quad \arbrebbc{\rond}{\lower2pt\hbox{$\;$}}{\raise2pt
\hbox{$\kern-2pt
\carre$}}{k+2}{}{}{k+1}{\raise5pt\hbox{$\scriptstyle\kern4pt
2n$}}{\textstyle\triang\kern-8pt\backslash}\raise3pt\hbox{$\scriptstyle k$}
{\raise10pt\hbox{$\scriptstyle\kern-8pt $}}\qquad
+\quad \arbrebbc{\rond}{\lower2pt\hbox{$\;$}}{\raise2pt
\hbox{$\kern-2pt
\carresec$}}{k+2}{}{}{k+1}{\raise5pt\hbox{$\scriptstyle\kern4pt
2n$}}{\textstyle\triang}\raise3pt\hbox{$\scriptstyle k$}
{\raise10pt\hbox{$\scriptstyle\kern-8pt $}}\qquad
\cr
\noalign{\smallskip}
&:=D(\star)+ D(\triang\kern-8pt\backslash\,)+ D(\carresec).\cr}
$$
In the above seven trunks of the two decompositions the symbols ``$\carre$,'' ``$\triang$,'' ``$\rond$'' without slash
are not empty and do not contain the one-child node.
Furthermore, note that
$D=D(\star)+ D(\triang\kern-8pt\backslash\,)+ D(\textstyle\rondsec)$ also holds.

\medskip
Pivoting the two subtrees on each of the nodes $k$, $(k+1)$, 
$(k+2)$ in $C(\star)$ (resp. in $D(\star)$) yields the same trunk, as in (2.1). We then say that {\it subtree pivoting is permitted} on those nodes. Furthermore, 
the three subtrees ``$\rond$,'' ``$\carre$'' and ``$\triang$'' in $C(\star)$ play a symmetric
role, while
in $D(\star)$ only ``$\rond$'' and ``$\triang$,'' on the one hand,
and ``$\carre$'' and ``$\triang$,'' on the other hand,
have a symmetric role.
Hence,
$$
C(\star)=2\,D(\star). 
$$

Now, remember that in each secant tree the one-child node lies at the rightmost position. 
This implies that {\it subtree pivoting on each of 
the ancestors of the one-child is not permitted}.
Thus, subtree pivoting in $C(\triangsec)$ and $D(\triangsec)$ is only permitted on nodes $k+2$ and $k+1$.
On the other hand, the two subtrees ``$\rond$'' and ```$\carre$'' play a symmetric role in $C(\triangsec)$, but not in $D(\triangsec)$.
Hence,
$$
C(\triangsec)=2\,D(\triangsec). 
$$
In $C(\carresec)$ and $D(\carresec)$
the two subtrees ``$\rond$'' and ``$\triang$'' play a symmetric
role. Moreover, subtree pivoting on nodes $k+2$ is permitted
in $C(\carresec)$, but not in $D(\carresec)$.
Hence,
$$\leqalignno{
C(\carresec)&=2\,D(\carresec),\cr 
\noalign{\hbox{so that}}
C-2\,D& =C(\rondsec).\cr}  
$$

\bigskip
\centerline{\bf 3. Proof that $(R\,1)$ holds}
The decomposition
$$
{\goth T}_{2n,m,k}
=\quad\arbrea{m}{\lower4.5pt\hbox{\kern-3pt$\carre\atop
m+1$}}{} {\lower2pt\hbox{$\bullet$}}\quad
+\quad[\quad\arbrea{m}{\rond}{}{\lower2pt\hbox{$\bullet$}},
\quad m+1]
$$
means that in each tree from
${\goth T}_{2n,m,k}$  the node $(m+1)$ is, or is not, the
sibling of the leaf~$m$. In the next decomposition the node~$m$
is, or is not, the parent of the leaf~$(m+1)$:
$$
{\goth T}_{2n,m+1,k}
=\quad\arbrea{m+1}{\rond}{m} {\lower2pt\hbox{$\bullet$}}\quad
+\quad[\quad\arbrea{m+1}{\rond}{}{\lower2pt\hbox{$\bullet$}}\quad,m].
$$
As already explained in Example 1, Section 2, we can write:
$$
\Deltaa_m{\kern-3pt}\,
{\goth T}_{2n,m,k}
=
{\goth T}_{2n,m+1,k}-{\goth T}_{2n,m,k}
=\arbrea{m+1}{\lower0pt\hbox{$\rond
$}}{m}{\lower2pt\hbox{$\bullet$}}\qquad
-\qquad\arbrea{m}{\lower4.5pt\hbox{\kern-5pt$\rond\atop
\ m+1$}}{} {\lower2pt\hbox{$\bullet$}}\quad,
$$

\goodbreak\noindent
so that
$$
\leqalignno{
\Deltaa_m{\kern-3pt}^2&\,
{\goth T}_{2n,m,k}\cr\noalign{\vskip-3pt}
&=
({\goth T}_{2n,m+2,k}-{\goth T}_{2n,m+1,k})
-({\goth T}_{2n,m+1,k}-{\goth T}_{2n,m,k})\cr
&=\quad\arbrea{m+2}{\raise1.5pt\hbox{\kern-3pt$\rond$}}{m+1} {\lower1.5pt\hbox{$\bullet$}}\qquad
-\qquad\arbrea{m+1}{\lower4.5pt\hbox{\kern-5pt$\rond\atop
\ m+2$}}{} {\lower2pt\hbox{$\bullet$}}\quad
-\qquad\arbrea{m+1}{\lower0pt\hbox{$\rond
$}}{m}{\lower2pt\hbox{$\bullet$}}\qquad
+\qquad\arbrea{m}{\lower4.5pt\hbox{\kern-5pt$\rond\atop
\ m+1$}}{} {\lower2pt\hbox{$\bullet$}}\quad\cr 
&:=\kern18pt A\kern30pt  -\kern30pt B
\kern18pt -\kern31pt C\kern27pt +\kern28pt D\quad
.\cr}$$
Depending on the mutual positions of nodes~$m$, $(m+1)$ and $(m+2)$ the further decompositions prevail, as again 
$k$ still remains attached to $(2n)$:
$$
\leqalignno{\noalign{\smallskip}
A&=\qquad\arbrebb{\lower2pt\hbox{$\;\bullet$}}{\carre}%
{\rond}{m+1}{m+2}{}{m}{}\quad
+\quad[\quad\arbrebb{\lower2pt\hbox{$\;\bullet$}}{\carre}%
{\rond}{m+1}{m+2}{}{}{}\quad,m]:=A_{1}+A_{2};\cr
\noalign{\medskip}
B&=\quad \arbrea{m+1}{\rond}{m}{\lower2pt
\hbox{$\bullet$}}\centerput(1,4){$\scriptstyle m+2$}\quad
+
\quad[\quad\arbrea{m+1}{\lower4.5pt\hbox{\kern-3pt$\rond\atop
m+2$}}{} {\lower2pt\hbox{$\bullet$}}\quad,m]
:=B_{1}+B_{2};\cr
\noalign{\medskip}
C&=\quad \arbrea{m+1}{\rond}{m}{\lower2pt
\hbox{$\bullet$}}\centerput(1,4){$\scriptstyle m+2$}\quad
+\qquad\arbrebb{\lower2pt\hbox{$\;\bullet$}}{\carre}%
{\rond}{m}{m+1}{}{}{m+2}\quad
+\quad[\quad\arbrebb{\lower2pt\hbox{$\;\bullet$}}{\carre}%
{\rond}{m}{m+1}{}{}{},m+2] :=C_{1}+C_{2}+C_{3};\cr
\noalign{\medskip}
D&=\quad\arbrebbj{}{m+2}{}%
{\bullet}{\carre}{\rond}{m}{m+1}\qquad
+\quad[\quad\arbrea{m}{\kern-5pt \lower3pt\hbox{$\rond\atop \ m+1$}}{}{\lower2pt
\hbox{$\bullet$}}\quad,m+2]:=D_{1}+D_{2}.\cr}
$$

In the above decompositions the subsets $B_{1}$ and $C_{1}$
are identical. Furthermore,
$A_{2}=C_{3}$ and $B_{2}=D_{2}$. Accordingly,
$\Deltaa_m{\kern-3pt}^2\,
{\goth T}_{2n+1,m,k}=A_{1}-2\,B_{1}-C_{2}+D_{1}$.
A further decomposition of those four terms, depending upon the occurrence of the one-child node, is to be worked out:
$$\leqalignno{
A_{1}&:=A_{1}(\star)+A_{1}(\carresec)+A_{1}(\rondsec);\cr
B_{1}&:=B_{1}(\star)+B_{1}(\rondsec);\cr
C_{2}&:=C_{2}(\star)+C_{2}(\carresec)+C_{2}(\rondsec);\cr
D_{1}&:=D_{1}(\star)+D_{1}(\carresec)+D_{1}(\rondsec).\cr
}
$$
Now, $B_{1}(\star)$ can be decomposed into
$$\leqalignno{
B_{1}(\star)=\quad \arbrea{m+1}{\rond}{m}{\lower2pt
\hbox{$\bullet$}}\centerput(1,4){$\scriptstyle m+2$}\quad
&=\quad \arbrea{m+1}{}{m}{\lower2pt
\hbox{$\bullet$}}\centerput(1,4){$\scriptstyle m+2$}\quad
+\quad
\arbrebbj{}{}{}{\textstyle\bullet}{\rond}{\carre}{m+1}{m+2}
\centerput(-2.5,2){$\scriptstyle m$}\qquad\cr
&=B_{1,1}(\star)+B_{1,2}(\star),\cr
}
$$
where ``$\carre$'' and ``$\rond$'' are supposed to be {\it non-empty} in $B_{1,2}(\star)$. 
By the Tree Calculus techniques developed in Section 2, we have
$A_{1}(\star)=2\,B_{1,2}(\star)$. 
On the other hand, $B_{1}(\rondsec)$ can also be written
$$
B_{1}(\rondsec)=
\quad
\arbrebbj{}{}{}{\textstyle\bullet}{\textstyle\rondsec}{\carre}{m+1}{m+2}
\centerput(-2.5,2){$\scriptstyle m$}\quad,
$$
so that
$A_{1}(\carresec)=2\,B_{1}(\rondsec)$ and $A_{1}(\rondsec)=B_{1}(\rondsec)$.
Furthermore, $C_{2}(\star)=D_{1}(\star)$, 
$C_{2}(\carresec)=2\,D_{1}(\carresec)$ and $C_{2}(\rondsec)=D_{1}(\rondsec)$. Altogether,
$A_{1}-2\,B_{1}-C_{2}+D_{1}=\bigl(A_{1}(\star)+A_{1}(\carresec)+A_{1}(\rondsec)\bigr)-2\bigl(B_{1,1}(\star)+B_{1,2}(\star)+B_{1}(\rondsec)\bigr)
-\bigl(C_{2}(\star)+C_{2}(\carresec)+C_{2}(\rondsec)\bigr)+\bigl(D_{1}(\star)+D_{1}(\carresec)+D_{1}(\rondsec)\bigr)=A_{1}(\rondsec)-2\,B_{1,1}(\star)-D_{1}(\carresec)$.

As $m$ is supposed to be at least equal to~2,  we can write
$$A_{1}(\rondsec)=\qquad
\quad\arbrebbc{\ \bullet m+2}{\carre}{}{m+1}{}{\kern -4pt\textstyle\rondsec\quad}{m}{}{\triang}
\raise3pt\hbox{}
\qquad,\quad
D_{1}(\carresec)=\qquad
\vbox{\vskip1.cm\offinterlineskip
\segment(0,0)\dir(0,1)\long{3}
\segment(0,3)\dir(-1,1)\long{4}
\segment(0,3)\dir(1,1)\long{12}
\segment(4,7)\dir(-1,1)\long{4}
\segment(8,11)\dir(-1,1)\long{4}
\centerput(-4,8){$\scriptstyle m$}
\centerput(-4,6){$\bullet$}
\centerput(0,11.5){$\rond$}
\centerput(12,15.5){$\carresec$}
\centerput(7.5,6){$\scriptstyle m+1$}
\centerput(12,10.5){$\scriptstyle m+2$}
\centerput(4,15.5){$\triang$}
}\qquad\qquad\quad,
$$
which shows that those two families are equal.
Thus,
$$\Deltaa_m{\kern-3pt}^2\,{\goth T}_{2n,m,k}
=-2\,B_{1,1}(\star)=-2\quad \arbrea{m+1}{}{m}{\lower2pt
\hbox{$\bullet$}}\centerput(1,5){$\scriptstyle m+2$}.$$ This expression is also equal to
${}-4\,{\goth T}_{2n-2,m,k-2}$, because in each tree~$t$ from $B_{1}$ the nodes $(m+1)$ and $(m+2)$ are both leaves. Remove them, as well as the two edges going out of~$m$, and subtract~2 from all the remaining nodes greater than $(m+2)$. The tree thereby derived belongs to 
${\goth T}_{2n-2,m,k-2}$.\qed

\def\arbrebb#1#2#3#4#5#6#7#8{\mathop{\hskip14pt
\vbox{\vskip1.2cm\offinterlineskip
\segment(0,0)\dir(0,1)\long{3}
\segment(0,3)\dir(1,1)\long{4}
\segment(0,3)\dir(-1,1)\long{8}
\segment(-4.5,7)\dir(1,1)\long{4}
\leftput(-10,11.5){$\scriptstyle#1$}
\leftput(2.5,7.6){\hbox{$\scriptstyle#2$}}
\leftput(-1.3,10.8){\hbox{$\scriptstyle#3$}}
\rightput(-4,4.5){\hbox{$\scriptstyle#4$}}
\rightput(-8,9){\hbox{$\scriptstyle#5$}}
\leftput(0.2,11){\hbox{$\scriptstyle#6$}}
\leftput(0.5,1.5){\hbox{$\scriptstyle#7$}}
\leftput(4,5){\hbox{$\scriptstyle#8$}}
}\hskip14pt}\nolimits}

\def\arbrebbj#1#2#3#4#5#6#7#8{\mathop{\hskip14pt
\vbox{\vskip1.2cm\offinterlineskip
\segment(0,0)\dir(0,1)\long{3}
\segment(0,3)\dir(-1,1)\long{4}
\segment(0,3)\dir(1,1)\long{8}
\segment(4.5,7)\dir(-1,1)\long{4}
\leftput(-10,11.5){$\scriptstyle#1$}
\leftput(7.5,9){\hbox{$\scriptstyle#2$}}
\leftput(-1.3,10.8){\hbox{$\scriptstyle#3$}}
\centerput(-4,6.5){\hbox{$\scriptstyle#4$}}
\rightput(10,11.5){\hbox{$\scriptstyle#5$}}
\leftput(-1.5,11.5){\hbox{$\scriptstyle#6$}}
\rightput(-4,4.5){\hbox{$\scriptstyle#7$}}
\leftput(4,5){\hbox{$\scriptstyle#8$}}
}\hskip14pt}\nolimits}

\bigskip
\centerline{\bf 4. Tree Calculus for proving that  $(R\,2)$ holds}

\medskip

\def\arbrebb#1#2#3#4#5#6#7#8{\mathop{\hskip14pt
\vbox{\vskip1.2cm\offinterlineskip
\segment(0,0)\dir(0,1)\long{3}
\segment(0,3)\dir(1,1)\long{4}
\segment(0,3)\dir(-1,1)\long{8}
\segment(-4.5,7)\dir(1,1)\long{4}
\leftput(-10,11.5){$\scriptstyle#1$}
\leftput(2.5,7.6){\hbox{$\scriptstyle#2$}}
\leftput(-1.3,10.8){\hbox{$\scriptstyle#3$}}
\rightput(-4,4.5){\hbox{$\scriptstyle#4$}}
\rightput(-8,9){\hbox{$\scriptstyle#5$}}
\leftput(0.2,11){\hbox{$\scriptstyle#6$}}
\leftput(0.5,1.5){\hbox{$\scriptstyle#7$}}
\leftput(4,5){\hbox{$\scriptstyle#8$}}
}\hskip14pt}\nolimits}

\def\arbrebbj#1#2#3#4#5#6#7#8{\mathop{\hskip14pt
\vbox{\vskip1.2cm\offinterlineskip
\segment(0,0)\dir(0,1)\long{3}
\segment(0,3)\dir(-1,1)\long{4}
\segment(0,3)\dir(1,1)\long{8}
\segment(4.5,7)\dir(-1,1)\long{4}
\leftput(-10,11.5){$\scriptstyle#1$}
\leftput(7.5,9){\hbox{$\scriptstyle#2$}}
\leftput(-1.3,10.8){\hbox{$\scriptstyle#3$}}
\centerput(-4,6.5){\hbox{$\scriptstyle#4$}}
\rightput(10,11.5){\hbox{$\scriptstyle#5$}}
\leftput(-1.5,11.5){\hbox{$\scriptstyle#6$}}
\rightput(-4,4.5){\hbox{$\scriptstyle#7$}}
\leftput(4,5){\hbox{$\scriptstyle#8$}}
}\hskip14pt}\nolimits}

With $n\ge 3$ and 
$2\le m\le k-1<k\le 2n-3$ we have:
$$
{\goth T}_{2n,m,k}
=[\arbrea{k+1}{\ 2n}k{\rond},m]+
[\arbrea{}{\ 2n}k{\rond}, m,k+1],
$$
meaning that each tree from ${\goth T}_{2n,m,k}$ has one of
the two forms: either $k+1$ is incident to~$k$, or not, and  the leaf~$m$ is the end of the minimal chain.

Using the same dichotomy,

$$
{\goth T}_{2n,m,k+1}
=[\quad\arbrebb{\rond}{\carre}{}{k+1}{}{2n}{k}{},m]
+[\arbrea{}{\ 2n}{k+1}{\rond}\quad,m, k].
$$

\goodbreak
As the second terms of the above two equations are in one-to-one correspondence by the transposition $(k,k+1)$ we have: 
$$\leqalignno{\noalign{\vskip5pt}
{\goth T}_{2n,m,k+1}-{\goth T}_{2n,m,k}
&=[\quad\arbrebb{\rond}{\carre}{}{k+1}{}{2n}{k}{},m]
-[\quad\arbrea{k+1}{\ 2n}k{\rond},m]
:=B-A.\cr
\noalign{\hbox{
In the same manner,}}
{\goth T}_{2n,m,k+2}-{\goth T}_{2n,m,k+1}
&=[\quad\arbrebb{\rond}{\carre}{}{k+2}{}{2n}{k+1}{}\ ,m]
-[\quad\arbrea{k+2}{\ 2n}{k+1}{\rond}\ ,m]
:=D-C.\cr}
$$ 
$$\leqalignno{\qquad
\Deltaa_k{\kern-3pt}^2\,{\goth T}_{2n,m,k}&=\bigl({\goth T}_{2n,m,k+2}-{\goth T}_{2n,m,k+1}\bigr)
-\bigl({\goth T}_{2n,m,k+1}-{\goth T}_{2n,m,k}\bigr)&\hbox{Thus,}\cr
&=D-C-B+A.\cr
}
$$

The further decompositions of the components of the previous sum depend on the mutual positions of the nodes~$k$, $(k+1)$, $(k+2)$;
$$
\leqalignno{
	D&=[\qquad\arbrebb{\rond}{\carre}{}{k+2}{}{2n}{k+1}{}\quad
,\ m]
=[\qquad\arbrebb{\rond}{\carre}{}{k+2}{}{2n}{k+1}{}\quad
,\ m,k\ ]+
[\quad\arbrebbc{\rond}{\carre}{}{k+2}{}{2n}{k+1}{}{\triang}
\raise3pt\hbox{$\scriptstyle k$}\quad,m]\cr
&\kern 3.2cm:=D_{1}+D_{2};\cr
\noalign{\smallskip}
C&=[\quad\arbrea{k+2}{\,2n}{k+1}{\rond}\quad
,\ m]
=[\quad\arbrea{k+2}{\,2n}{k+1}{\rond}\quad
,\ m,k\ ]\quad
+\quad[\quad\arbrebb{\rond}{\lower2.5pt\hbox{$\scriptstyle\
\carre$}}{}{k+1}{k+2}{2n}{k}{}\quad,\ m]\cr
&\kern 2.9cm:=C_{1}+C_{2};\cr
\noalign{\medskip\smallskip}
B&=[\quad\arbrebb{\rond}{\carre}{}{k+1}{}{2n}{k}{},\ m]
=[\quad\arbrebb{\rond}{\carre}{}{k+1}{}{2n}{k}{},\ m,k+2]
\quad
+[\quad\arbrebb{\lower2pt\hbox{$\;$}}{\carre}%
{}{k+1}{k+2}{2n}{k}{},m]\quad
\cr
&\kern 1.5cm{}
+\quad
[\arbrebbc{\rond}{\lower2pt\hbox{$\;$}}{\raise2pt
\hbox{$\kern-2pt
\carre$}}{k+2}{}{}{k+1}{\raise5pt\hbox{$\scriptstyle\kern4pt
2n$}}{\triang}\raise3pt\hbox{$\scriptstyle k$}
{\raise10pt\hbox{$\scriptstyle\kern-8pt $}}\quad,m]
+\quad[\quad \arbrebb{\carre}{\lower3pt\hbox{$\ $}}{}{k+1}{}{2n}{k}{k+2}\quad,m]
+\quad[\quad
\lower9pt\hbox{$\arbrebbz{\triang}{\quad k+2}{}{k+1}{2n}{\rond}{k}{}{\carre}$}\quad,m]\cr
&\kern 2.6cm:=B_{1}+B_{2} +B_{3}+B_{4}+B_{5};\cr
\noalign{\vskip5pt}
A&=[\quad\arbrea{k+1}{\,2n}{k}{\rond}
,\ m]
=[\quad\arbrea{k+1}{\,2n}{k}{\rond}
,\ m,k+2\ ]+[\qquad\arbrebb{\rond}{\lower2.5pt\hbox{$\scriptstyle\
$}}{\carre}{k+1}{k+2}{}{k}{2n}\quad,m]\cr
&\kern 2.5cm:=A_{1}+A_{2},\cr
}
$$
where ``$\rond$'' and ``$\carre$'' cannot be  both empty in~$B_{3}$ and $B_{5}$.

\goodbreak
Obviously, $D_{1}=B_{1}$ and~$C_{1}=A_{1}$, so that the sum
$D-C-B+A$ may be written 
$(B_{2}-B_{4})-(C_{2}-A_{2})+(D_{2}-B_{3}-B_{5})-2\,B_{2}$, showing that the components fall in three categories: (1) $B_{2}$ and $B_{4}$ having one subtree to determine; (2) $C_{2}$ and $A_{2}$ having two such subtrees;
(3) $D_{2}$, $B_{3}$ and $B_{5}$ having three of them. 

(1) With the same notation as in (2.2) write: $B_{2}=B_{2}(\star)+B_{2}(\carresec)$ and $B_{4}=B_{4}(\star)+B_{4}(\carresec)$. First, $B_{2}(\star)=B_{4}(\star)$, but
$B_{2}(\carresec)=2\,B_{4}(\carresec)$. Hence,
$$
B_{2}-B_{4}=B_{4}(\carresec)=
[\quad\arbrebbj{}{}{}%
{}{\carresec}{2n}{k+2}{k+1}
\centerput(-4,2){$\scriptstyle k$}\quad,m].
$$

(2) Next, $C_{2}=C_{2}(\star)+C_{2}(\carresec)+C_{2}(\rondsec)$ and
$A_{2}=A_{2}(\star)+A_{2}(\carresec)+A_{2}(\rondsec)$. As before, $C_{2}(\star)=A_{2}(\star)$ and also 
$C_{2}(\rondsec)=A_{2}(\rondsec)$. However, $C_{2}(\carresec)=2\,A_{2}(\carresec)$. Altogether,
$$
C_{2}-A_{2}=A_{2}(\carresec)=
[\quad\arbrebbj{}{}{}%
{}{\carresec}{\rond}{2n}{k+1}
\centerput(-4,2){$\scriptstyle k$}
\centerput(-7.5,8.5){$\scriptstyle k+2$}\quad,m].
$$

(3) The calculation of the third sum requires the following decompositions:
$$\leqalignno{
D_{2}
&:=D_{2}(\star)+D_{2}(\triang\kern-8pt\backslash\,)+D_{2}(\carresec)+D_{2}(\rondsec);\cr
B_{3}
&:=B_{3}(\star)+ B_{3}(\triang\kern-8pt\backslash\,)+ B_{3}(\carresec);\cr
B_{5}
&:=B_{5}(\star)+ B_{5}(\triang\kern-8pt\backslash\,)+ B_{5}(\carresec);\cr}
$$
By the Tree Calculus techniques developed in Section 2, especially Example 2, we have
$$
\leqalignno{
	D_{2}(\star) & =2\,B_{3}(\star)=2\,B_{5}(\star);\cr
	D_{2}(\triang\kern-8pt\backslash\,) &=2\,B_{3}(\triang\kern-8pt\backslash\,)
	=2\,B_{3}(\triang\kern-8pt\backslash\,);\cr
D_{2}(\carresec)&=2\,B_{3}(\carresec)=2\,B_{3}(\carresec).\cr
}$$
Hence,

$$D_{2}-B_{3}-B_{5}
=D_{2}(\rondsec)
=[\qquad
\vbox{\vskip1.cm\offinterlineskip
\segment(0,0)\dir(0,1)\long{3}
\segment(0,3)\dir(-1,1)\long{4}
\segment(0,3)\dir(1,1)\long{12}
\segment(4,7)\dir(-1,1)\long{4}
\segment(8,11)\dir(-1,1)\long{4}
\centerput(-4,7.5){$\triang$}
\centerput(0,11.5){$\carre$}
\centerput(2,2.5){$\scriptstyle k$}
\centerput(12,16){$\rondsec$}
\centerput(7.5,6.5){$\scriptstyle k+1$}
\centerput(12,10.5){$\scriptstyle k+2$}
\centerput(4.5,16.5){$\scriptstyle 2n$}
}\qquad\quad,m].\hskip0.2cm 
$$

Now
$$
A_{2}(\carresec)=
[\quad\arbrebbj{}{}{}%
{}{\carresec}{\rond}{2n}{k+1}
\centerput(-4,2){$\scriptstyle k$}
\centerput(-7.5,8.5){$\scriptstyle k+2$}\quad,m]
=
[\quad\arbrebbj{}{}{}%
{}{\carresec}{}{2n}{k+1}
\centerput(-4,2){$\scriptstyle k$}
\centerput(-7.5,8.5){$\scriptstyle k+2$}\quad,m]
+
[\qquad
\vbox{\vskip1.cm\offinterlineskip
\segment(0,0)\dir(0,1)\long{3}
\segment(0,3)\dir(-1,1)\long{4}
\segment(0,3)\dir(1,1)\long{8}
\segment(4,7)\dir(-1,1)\long{8}
\segment(0,11)\dir(1,1)\long{4}
\centerput(-5,4.5){$\scriptstyle 2n$}
\centerput(-4,15.5){$\rond$}
\centerput(8,11.5){$\carresec$}
\centerput(7.5,6){$\scriptstyle k+1$}
\centerput(2,2){$\scriptstyle k$}
\centerput(-3.5,10){$\scriptstyle k+2$}
\centerput(4,15.5){$\triang$}
}\qquad\quad,m],
$$
where the two subtrees ``$\rond$'' and ``$\triang$'' are non-empty.
The last identity can be rewritten: $A_{2}(\carresec)=B_{4}(\carresec)+ D_{2}(\rondsec)$.
Hence,
$\Deltaa_k{\kern-3pt}^2\,{\goth T}_{2n,m,k}=D-C-B+A=
(B_{2}-B_{4})-(C_{2}-A_{2})+(D_{2}-B_{3}-B_{5})-2\,B_{2}=
B_{4}(\carresec)-A_{2}(\carresec)+D_{2}(\rondsec)-2\,B_{2}=-2\,B_{2}=-2\,[\quad\arbrebb{\lower2pt\hbox{$\;$}}{\carre}%
{}{k+1}{k+2}{2n}{k}{},m].
$
The last term is also equal to
$-4\,{\goth T}_{2n-2,m,k}$, because in each tree~$t$ from $B_{2}$ the nodes $(k+2)$ and $(2n)$ are both leaves. Remove them, as well as the two edges going out of~$(k+1)$, change $(k+1)$ into $(2n-2)$ and subtract~2 from all the remaining nodes greater than $(k+2)$. The tree thereby derived belongs to 
${\goth T}_{2n-2,m,k}$.\qed

\bigskip
\centerline{\bf 5. Proofs of Theorem 1.2, Proposition 1.3
and Theorem 1.4}

\medskip
{\it Proof of Theorem $1.2$. 

First top row}.\quad  The trees $t$ from ${\goth T}_{2n}$ such that $\eoc (t)=2$ contain the edge $2\rightarrow 1$. Remove it and change each remaining node label~$j$ by $j-2$. We get a tree $t'$ from
${\goth T}_{2n-2}$ such that $\pom (t')=k-2$. The second equality is a consequence of (1.4).\qed

\medskip
{\it Second top row}.\quad 
As $k\ge 4$, the subtrees ``$\carre$'' and ``$\rond$'' below are non-empty, so that
$$\leqalignno{\noalign{\vskip-14pt}
f_{2n}(3,k)&=\qquad\qquad
\vbox{\vskip1.4cm\offinterlineskip
\segment(0,0)\dir(-1,1)\long{10}
\segment(0,0)\dir(1,1)\long{5}
\segment(-5,5)\dir(1,1)\long{5}
\centerput(0,-3){$\scriptstyle1$}
\centerput(0,7){\hbox{$\scriptstyle 3$}}
\centerput(-5,2){$\scriptstyle 2$}
\centerput(0,9){$\scriptstyle \bullet$}
\centerput(-10,10.5){$\rond$}
\centerput(5,5.5){$\carre$}
}
\qquad=
\qquad\qquad
\vbox{\vskip1.4cm\offinterlineskip
\segment(0,0)\dir(-1,1)\long{10}
\segment(0,0)\dir(1,1)\long{5}
\segment(-5,5)\dir(1,1)\long{5}
\centerput(0,-3){$\scriptstyle1$}
\centerput(0,7){\hbox{$\scriptstyle 3$}}
\centerput(-5,2){$\scriptstyle 2$}
\centerput(0,9){$\scriptstyle \bullet$}
\centerput(-10,10.5){$\rond$}
\centerput(5,5.5){$\carresec$}
}
\qquad+\qquad\quad
\vbox{\vskip1.4cm\offinterlineskip
\segment(0,0)\dir(-1,1)\long{5}
\segment(0,0)\dir(1,1)\long{10}
\segment(5,5)\dir(-1,1)\long{5}
\centerput(0,-3){$\scriptstyle1$}
\centerput(0,7){\hbox{$\scriptstyle 3$}}
\centerput(5,2){$\scriptstyle 2$}
\centerput(0,9){$\scriptstyle \bullet$}
\centerput(-5,5.5){$\carre$}
\centerput(10,10.5){$\rondsec$}
}\qquad\quad\cr
\noalign{\bigskip}
&:=\kern 1cm A\kern.7cm{}=\kern 1cm A(\carresec\,)
\kern.3cm+\kern.7cm A(\rondsec);\cr
f_{2n}(2,k)&=\qquad\vbox{\vskip1.4cm\offinterlineskip
\segment(0,0)\dir(-1,1)\long{5}
\segment(0,0)\dir(1,1)\long{10}
\segment(5,5)\dir(-1,1)\long{5}
\centerput(0,-3){$\scriptstyle1$}
\leftput(6,3){\hbox{$\scriptstyle 3$}}
\rightput(-6,5){$\scriptstyle 2$}
\centerput(-5,4){$\scriptstyle \bullet$}
\centerput(0,10.5){$\rond$}
\centerput(10,10.5){$\carre$}
}\qquad\qquad:=B=\qquad\vbox{\vskip1.4cm\offinterlineskip
\segment(0,0)\dir(-1,1)\long{5}
\segment(0,0)\dir(1,1)\long{10}
\segment(5,5)\dir(-1,1)\long{5}
\centerput(0,-3){$\scriptstyle1$}
\leftput(6,3){\hbox{$\scriptstyle 3$}}
\rightput(-6,5){$\scriptstyle 2$}
\centerput(-5,4){$\scriptstyle \bullet$}
\centerput(0,10.5){$\rond$}
\centerput(10,10.5){$\carresec$}
}\qquad\qquad=\quad B(\carresec),\cr}$$
since $B(\star)$ is empty.
Therefore, $A(\carresec\,)=2\,B(\carresec\,)$, $A(\rondsec)=
B(\carresec\,)$ and then $f_{2,n}(3,k)=3\,f_{2n}(2,k)$.\qed

\bigskip
{\it Rightmost column}.\quad
The trees $t$ from ${\goth T}_{2n}$ such that $\eoc(t)=m$ and
$\pom(t)=2n-1$ contain the {\it rightmost} path $2n \rightarrow
(2n-1)\rightarrow{}$. Remove it. What is left is a tree $t'$ from ${\goth T}_{2n-2}$ such that $\eoc(t')=m$. The second equality is a consequence of (1.4).\cqfd

\goodbreak
\medskip
{\it Next to rightmost column}.\ 
If $t\in {\goth T}_{2n}$ and $\pom(t)=2n-1$, then~$t$~contains
the rightmost path $2n\rightarrow (2n-1)\rightarrow$, as already noted. On the other hand, the node with label
$(2n-2)$ is necessarily a leaf. The transposition $(2n-2,2n-1)$
transforms~$t$ into a tree~$t_1$ such that $\pom(t_1)=2n-2$.
Also, removing the path $ 2n\rightarrow (2n-1)$
and rooting either the subtree
${}^{2n-1}{\nwarrow\!\nearrow}^{2n}$, or the subtree
${}^{2n}{\nwarrow\!\nearrow}^{2n-1}$, onto the leaf labeled
$(2n-2)$ gives rise to two trees $t_2$, $t_3$ such that
$\pom(t_2)=\pom(t_3)=2n-2$.\qed

\bigskip
{\it Proof of Proposition $1.3$}.\quad
For (1.8) and (1.9) we only have to reproduce the proofs made in Theorem~1.2 for the first and second top rows, the parent of the maximum node playing no role.

To obtain $(R\,4)$ simply write $(R\,2)$ for $m=2$ and rewrite it using the first identity of Theorem~1.2 dealing with the ``first top row.'' Next, $(R\,3)$ is deduced from $(R\,4)$ by using identity (1.4).\qed

\bigskip
{\it Proof of Theorem $1.4$}.\quad
By induction, the row and column sums $f_{2n}(m,\brullet)$ and $f_{2n}(\brullet,k)$ of the matrices $M_{2n}$ can be calculated by means of relations (1.8), (1.9) and, either $(R\,3)$, or $(R\,4)$. The first and second top rows (resp. rightmost and next to rightmost columns) of the matrix~$M_{2n}$ are known by Theorem~1.2. It then suffices to apply, either rule $(R\,1)$, from top to bottom, or rule $(R\,2)$ from right to left to obtain the remaining entries of the upper triangles of~$M_{2n}$.\qed

\bigskip
\centerline{\bf 6. The lower triangles}

\medskip
Observe that for $2n\ge 6$ the non-zero entries of the first row, first column and last column of $M_{2n}$ are identical and they differ from the entries in the bottom row. For instance, the sequence $5,15,21,15,5$ in $M_8$ appears three times, but the bottom row reads: $16,16,14,10,5$. 

By Theorem 1.2 we already know that $f_{2n}(2,k)=
f_{2n-2}(k-1,\brullet)=f_{2n}(k-1,2n-1)$ $(3\le k\le 2n-1)$. We also have: $f_{2n}(m,1)=f_{2n-2}(m-1,\brullet)$ $(3\le m\le 2n-1)$, by using this argument: each tree $t$ from
${\goth T}_{2n}$ satisfying $\pom(t)=1$ has its leaf node~$2n$ incident to root~1. Just remove the edge ${2n}\rightarrow 1$. Change each remaining label~$j$ to $j-1$. We get a secant tree~$t'$ belonging to ${\goth T}_{2n-2}$. The mapping $t\mapsto t'$ is bijective; moreover, $\eoc(t')=\eoc(t)-1$.
As a summary,
$$
f_{2n}(2,k)=f_{2n}(k-1,2n-1)=f_{2n}(k,1)
\quad(3\le k\le 2n-1).\leqno(6.1)
$$

\def\ent{\mathop{\rm ent}\nolimits}
\def\Ent{\mathop{\rm Ent}\nolimits}

Introduce a further statistic~``$\ent$'' (short hand for ``Entringer'') on each tree~$t$ from ${\goth T}_{n}$ (even for {\it tangent} trees) as follows: $\ent(t)$ is the label of the {\it rightmost} node of~$t$. The
distribution of~``ent'' is well-known, mostly associated with the model of the alternating permutations and traditionally called
{\it Entringer} distribution (see, e.g., [En63], [Po82], [Po87], [GHZ10]). We also know how to calculate the
generating function of that distribution and  build up the
{\it Entringer triangle} $({\rm Ent}_n(j))$ $(1\le j\le
n-1;\, n\ge 2)$, as done in Table~6.1.

For instance, the leftmost 61 on the row $n=7$ is the sum of all the entries in the
previous row (including an entry~0 on the right!); the
second~61 is the sum of the leftmost five entries; 56, the sum of
the leftmost four entries; 46, the sum of the leftmost three 
entries; 32, the sum of the leftmost two entries and 16 is equal
to the leftmost entry. Let
$\Ent_n(j)=\#\{t\in {\goth T}_n:\ent(t)=j\}$.

$$
\matrix{j=&\kern8.5pt\vrule height 8pt width.4pt depth 4pt
&1&2&3&4&5&6\cr
\noalign{\hrule}
n=&2\ \vrule height 8pt width.4pt depth 4pt&1\cr
&3\ \vrule height 8pt width.4pt depth 4pt&1&1\cr
&4\ \vrule height 8pt width.4pt depth 4pt&2&2&1\cr
&5\ \vrule height 8pt width.4pt depth 4pt&5&5&4&2\cr
&6\ \vrule height 8pt width.4pt depth 4pt&16&16&14&10&5\cr
&7\ \vrule height 8pt width.4pt depth
4pt&61&61&56&46&32&16\cr }
$$

\centerline{Table 6.1. The Entringer distribution}

\proclaim Proposition 6.1. For $2\le k\le 2n-2$ we have:
$$f_{2n}(2n,k)=\Ent_{2n-2}(k-1).\leqno(6.2)
$$

\proof Note that in a tree $t\in {\goth
T}_{2n}$ such that $\eoc(t)=2n$ and $\pom(t)=k$ the leaf
labeled~$2n$ is the unique son of the node labeled~$k$, which
is also the rightmost node. Let
$(1=a_1)\rightarrow a_2\rightarrow a_3\rightarrow\cdots
\rightarrow (a_{j-1}=k)\rightarrow (a_j=2n)$ be the minimal
chain of~$t$. Form a new tree~$t'$ by means of the following
changes:

(i) delete the path ${}\rightarrow a_{j-1}\rightarrow 2n$;

(ii) for $i=1,2,\ldots ,j-2$ replace each node label $a_i$ of the
minimal chain by
$a_{i+1}-1$;

(iii) replace each other node label $b$ by $b-1$. 

\noindent
The label of the rightmost node of~$t'$ is then equal to the
$\pom(t)-1=k-1$. The mapping $t\mapsto t'$ is obvious
bijective.\cqfd

\medskip
Finally, the upper diagonal $\{f_{2n}(k+1,k):(1\le k\le 2n-1)\}$ of the lower triangle in each matrix $M_{2n}$ can also be fully evaluated. First, in an obvious manner,
$$
f_{2n}(2,1)=f_{2n}(2n,2n-1)=0. \leqno(6.3)
$$
Second, the identities
$$
f_{2n}(3,2)=2\,f_{2n}(3,1)=f_{2n}(2n-1,2n-2)=2\,f_{2n}(2n,2n-2)\leqno(6.4)
$$
(also equal to $f_{2n-4}(\brullet,\brullet)=E_{2n-4}$) can be proved as follows. To each tree $t\in {\goth T}_{2n}$ such that $\eoc(t)=3$, $\pom(t)=1$ there correspond {\it two} trees, whose ``eoc'' and ``pom'' are equal to~2 and~3, respectively, as illustrated by the diagram:
$$f_{2n}(3,1)=\qquad
\vbox{\vskip1.3cm\offinterlineskip
\segment(0,0)\dir(-1,1)\long{5}
\segment(0,0)\dir(1,1)\long{10}
\segment(5,5)\dir(-1,1)\long{5}
\centerput(0,-3){$\scriptstyle1$}
\centerput(0,7){\hbox{$\scriptstyle 3$}}
\centerput(-5,2){$\scriptstyle 2n$}
\centerput(5,2){$\scriptstyle 2$}
\centerput(0,9){$\scriptstyle \bullet$}
\centerput(10,10.5){$\carre$}
}\qquad\quad\mapsto
\qquad\qquad
\vbox{\vskip1.3cm\offinterlineskip
\segment(0,0)\dir(-1,1)\long{10}
\segment(0,0)\dir(1,1)\long{5}
\segment(-5,5)\dir(1,1)\long{5}
\centerput(0,-3){$\scriptstyle1$}
\centerput(-10.5,7){\hbox{$\scriptstyle 2n$}}
\centerput(0,7){\hbox{$\scriptstyle 3$}}
\centerput(-5,2){$\scriptstyle 2$}
\centerput(0,9){$\scriptstyle \bullet$}
\centerput(5,5.5){$\carre$}
}\qquad=f_{2n}(3,2).\qquad
$$
Likewise, each tree $t\in {\goth T}_{2n}$ such that 
$\eoc(t)=2n$ and $\pom(t)=2n-2$ necessarily has
its rightmost four nodes equal to $(2n-1), b,2n,(2n-2)$ in that order, with~$b$ being a node less than $(2n-2)$. In particular, the latter node is its rightmost (one-child) node.
To such a tree there correspond two trees, whose ``eoc'' and ``pom'' are equal to $(2n-1)$, $(2n-2)$, respectively, as illustrated by the next diagram. In particular, the node~$b$ becomes the rightmost node of those two such trees.
$$
\vbox{\vskip1.7cm\offinterlineskip
\segment(10,10)\dir(-1,1)\long{5}
\segment(0,0)\dir(1,1)\long{10}
\segment(5,5)\dir(-1,1)\long{5}
%
\centerput(-4,7.5){$\scriptstyle 2n-1$}
\centerput(2,14){$\scriptstyle 2n$}
\centerput(5,2){$\scriptstyle b$}
\centerput(5,14){$\scriptstyle \bullet$}
\centerput(14,10){$\scriptstyle 2n-2$}
}\qquad\qquad\qquad
\mapsto\qquad\qquad
\vbox{\vskip1.3cm\offinterlineskip
\segment(0,10)\dir(1,1)\long{5}
\segment(0,0)\dir(1,1)\long{5}
\segment(5,5)\dir(-1,1)\long{10}
%
\centerput(-4,8){$\scriptstyle 2n-2$}
\centerput(9.5,14){$\scriptstyle 2n-1$}
\centerput(5,2){$\scriptstyle b$}
\centerput(5,14){$\scriptstyle \bullet$}
\centerput(-7,14){$\scriptstyle 2n$}
}\qquad.\qquad\qed
$$

\proclaim Proposition 6.2. We have the crossing equalities:
$$
f_{2n}(k-1,k)+
f_{2n}(k+1,k)\!=\!
f_{2n}(k,k-1)+
f_{2n}(k,k+1),\leqno(6.5)
$$
for $3\le k\le 2n-2$.

The involved entries are located on the
four bullets drawn in the following diagramme.

$$
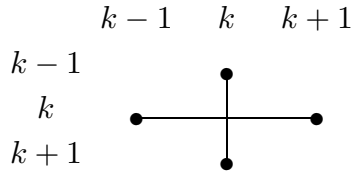
\vbox{\vskip1.3cm\offinterlineskip
\centerput(0,10){$k-1$}
\centerput(10,10){$k$}
\centerput(20,10){$k+1$}
\centerput(10,4){$\bullet$}
\centerput(20,-1){$\bullet$}
\segment(0,0)\dir(1,0)\long{20}
\centerput(0,-1){$\bullet$}
\centerput(10,-6){$\bullet$}
\segment(10,5)\dir(0,-1)\long{10}
\centerput(-10,5){$k-1$}
\centerput(-10,0){$k$}
\centerput(-10,-5){$k+1$}
\vskip.7cm
}\hskip1cm
$$

\proof
Let $i$, $j$ be two different integers from the set
$\{(k-1), k,(k+1)\}$. Say that $i$ and~$j$ are {\it connected} in a tree~$t$, if the tree contains the edge $i$---$j$, or 
if $i$ and $j$ are brothers and of them is the end of the minimal chain of~$t$. 
Each of the four ingredients of the previous identity 
is now decomposed into five terms, depending on whether the nodes $(k-1)$, $k$, $(k+1)$ are connected or not, namely: no connectedness; only $k, (k+1)$ connected; $(k-1), k$ connected; $(k-1),(k+1)$ connected; all connected. Thus,
$$
\leqalignno{
f_{2n}(k-1,k)&=
[\quad\arbrea{k-1}{\lower0pt\hbox{$\rond$}}{}{\lower2pt
\hbox{$\bullet$}}\ ,\quad
\arbrea{}{\lower0pt\hbox{$\carre$}}{k}{\raise2pt
\hbox{$\scriptstyle 2n$}}
,k+1]
+[\quad\arbrea{k-1}{\lower0pt\hbox{$\rond$}}{}{\lower2pt
\hbox{$\bullet$}}\ ,\qquad
\arbrea{2n}{\raise1pt\hbox{$\carre$}}{k}{\lower2pt
\hbox{}}\kern-5pt\raise12pt\hbox{$\scriptstyle
k+1$}\quad]\cr 
\noalign{\medskip}
&\kern-10pt{}+
[\quad
\arbrebb{\lower0pt\hbox{$\rond$}}{}{\raise4pt\hbox{$\scriptstyle
2n$}}{k}{}{}{}{\raise5pt\hbox{$\scriptstyle\bullet\  k-1$}}\ ,\ k+1\ ]
+[\quad\arbrea{k-1}{\lower4pt\hbox{$\rond\ \atop k+1$}}{}{\lower2pt
\hbox{$\bullet$}}\ ,\quad
\arbrea{}{\lower0pt\hbox{$\carre$}}{k}{\raise2pt
\hbox{$\scriptstyle 2n$}}]
+\qquad
\arbrebb{}{}{\raise4pt\hbox{$\scriptstyle
2n$}}{k}{\hbox{$\rond\atop k+1$}\hskip-5pt}{}{}{\raise5pt\hbox{$\scriptstyle\bullet\ k-1$}}
\cr
&:=A_1+A_2+A_3+A_4+A_5;\cr
f_{2n}(k+1,k)&=
[\ \arbrea{k+1}{\lower0pt\hbox{$\rond$}}{}{\lower2pt
\hbox{$\bullet$}}\ ,k-1,\arbrea{2n}{\lower0pt\hbox{$\rond$}}{k}{}]\quad
+\quad[\
\arbrea{k+1}{\lower-5pt\hbox{$\kern-15pt\textstyle\  \atop
2n$}}{k}{\lower2pt
\hbox{$\bullet$}}\
,k-1]\cr
\noalign{\medskip}
&\kern-1.5cm{}
+[\ \arbrea{k+1}{\lower0pt\hbox{$\rond$}}{}{\lower2pt
\hbox{$\bullet$}}\
,\qquad\quad \arbrebb{\lower2pt\hbox{}}{\rond}%
{}{k}{2n}{\kern-3pt \carre}{k-1}{}\quad]
+[\
\arbrea{k+1}{\lower0pt\hbox{$\carre$}}{k-1}{\lower2pt
\hbox{$\bullet$}}\
,\qquad
\arbrea{2n}{\lower0pt\hbox{$\rond$}}{k}{\lower2pt
\hbox{}}]\quad
+\quad\arbrebb{\lower2pt\hbox{$\;\bullet$}}{\rond}%
{}{k}{k+1}{\kern-3pt{\textstyle} 2n}{k-1}{}
\cr
\noalign{\smallskip}
&:=B_1+B_2+B_3+B_4+B_5;\cr
\noalign{\smallskip}
f_{2n}(k,k-1)&=
[\ \arbrea{k}{\lower0pt\hbox{$\rond$}}{}{\lower2pt
\hbox{$\bullet$}}\ ,k+1,\arbrea{2n}{\lower0pt\hbox{$\rond$}}{k-1}{}\quad]\quad
+[\ \arbrea{k}{\lower0pt\hbox{$\carre$}}{\raise5pt\hbox{$\kern5pt
\scriptstyle k+1$}}{\lower2pt
\hbox{$\bullet$}}\quad 
,\quad
\arbrea{2n}{\lower0pt\hbox{$\rond$}}{k-1}{\lower2pt
\hbox{}}\quad ]\cr
\noalign{\medskip}
&{}
+[\
\arbrea{k}{\lower-5pt\hbox{$\kern-15pt\textstyle\  \atop
2n$}}{k-1}{\lower2pt
\hbox{$\bullet$}}\
,k+1]
+[\ \arbrea{2n}{\lower4pt\hbox{$\carre\ \atop k+1$}}{k-1}{}\qquad 
,\ \arbrea{k}{\lower0pt\hbox{$\rond$}}{}{\lower2pt
\hbox{$\bullet$}}\ ]
\cr
\noalign{\smallskip}
&:=C_1+C_2+C_3+C_4;\cr
\noalign{\medskip}
f_{2n}(k,k+1)&=
[\quad\arbrea{k}{\lower0pt\hbox{$\rond$}}{}{\lower2pt
\hbox{$\bullet$}}\ ,\quad
\arbrea{}{\lower0pt\hbox{$\carre$}}{k+1}{\raise2pt
\hbox{$\scriptstyle 2n$}}
,\ k-1]
+[\quad
\arbrebb{\lower0pt\hbox{$\rond$}}{}{\raise4pt\hbox{$\scriptstyle
2n$}}{k+1}{}{}{}{\raise5pt\hbox{$\scriptstyle\bullet\  k$}}\ ,\ k-1\ ]
\cr 
\noalign{\medskip}
&\kern-25pt{}+
[\quad\arbrea{k}{\lower0pt\hbox{$\rond$}}{k-1}{\lower2pt
\hbox{$\bullet$}}\ ,\quad
\arbrea{2n}{\raise1pt\hbox{$\carre$}}{k+1}{\lower2pt
\hbox{}}\quad]
+[\quad\arbrebb{\lower0pt\hbox{$\carre$}}{}{\raise4pt\hbox{$\scriptstyle
2n$}}{k+1}{}{}{k-1}{\raise5pt\hbox{$\rond$}}\ ,\quad
\arbrea{k}{\lower0pt\hbox{$\carre$}}{}{\lower2pt
\hbox{$ \bullet$}}]
+\qquad
\arbrebb{}{}{\raise4pt\hbox{$\scriptstyle
2n$}}{k+1}{\raise5pt\hbox{$\rond$}}{}{k-1}{\raise5pt\hbox{$\scriptstyle\bullet\ k$}}
\cr
&:=D_1+D_2+D_3+D_4+D_5.\cr
}
$$

Now, the following identities hold: $A_1= C_1$, 
$A_2= C_4$, $A_3= D_{2}$,
$A_4= C_{2}$,  $B_1= D_1$,
$B_3= D_{4}$, $B_4= D_{3}$, so that
$\sum_{i}(A_{i}+B_{i})-\sum_{i}(C_{i}+D_{i})
=(B_{5}-D_{5}) -(C_{3}-A_{5}-B_{2})$. 

As before, we may write $B_{5}=B_{5}(\star)+B_{5}(\rondsec)$,
$D_{5}= D_{5}(\star)+ D_{5}(\rondsec)$. As 
$B_{5}(\star)=D_{5}(\star)$ and 
$B_{5}(\rondsec)=2\,D_{5}(\rondsec)$, we get:

\bigskip
$$
B_{5}-D_{5}=D_{5}(\rondsec)
=\arbrebbj{}{}{}%
{\bullet}{\textstyle\rondsec}{2n}{k}{k+1}
\centerput(-2,2){$\scriptstyle k-1$}\qquad
=\qquad\vbox{\vskip1.cm\offinterlineskip
\segment(0,0)\dir(0,1)\long{3}
\segment(0,3)\dir(-1,1)\long{4}
\segment(0,3)\dir(1,1)\long{12}
\segment(4,7)\dir(-1,1)\long{4}
\segment(8,11)\dir(-1,1)\long{4}
\centerput(-4,7.5){$\carre$}
\centerput(0,11.5){$\scriptstyle k$}
\centerput(0,9.5){$\scriptstyle\bullet$}
\centerput(12,16){$\rondsec$}
\centerput(8,6.5){$\scriptstyle k-1$}
\centerput(12,10.5){$\scriptstyle k+1$}
\centerput(4.5,15.5){$\scriptstyle 2n$}
}\qquad\quad,\leqno(6.6)
$$
as $k$ is supposed to be greater than~3.

\goodbreak
Next, 
$$\leqalignno{
C_{3}-B_{2}&=
[\
\arbrea{k}{\lower-5pt\hbox{$\kern-15pt\textstyle\  \atop
2n$}}{k-1}{\lower2pt
\hbox{$\bullet$}}\
,k+1]
-\ [\
\arbrea{k+1}{\lower-5pt\hbox{$\kern-15pt\textstyle\  \atop
2n$}}{k}{\lower2pt
\hbox{$\bullet$}}\
,k-1]\cr
\noalign{\medskip}
 &\kern-15pt=
\quad\arbrebb{\lower2pt\hbox{$\;\bullet$}}{\rond}%
{}{k-1}{k}{\kern-3pt{\textstyle} 2n}{}{k+1}
-\quad[\quad\arbrebb{\lower2pt\hbox{$\;\bullet$}}{\rond}%
{}{k-1}{k}{\kern-3pt{\textstyle} 2n}{}{},k+1]
-[\quad\arbrebb{\lower2pt\hbox{$\;\bullet$}}{\rond}%
{}{k}{k+1}{\kern-3pt{\textstyle} 2n}{}{},k-1]
=\quad\arbrebb{\lower2pt\hbox{$\;\bullet$}}{\rond}%
{}{k-1}{k}{\kern-3pt{\textstyle} 2n}{}{k+1}\quad;\cr
\noalign{\bigskip}
C_{3}-B_{2}&-A_{5}
=
\quad\arbrebb{\lower2pt\hbox{$\;\bullet$}}{\rond}%
{}{k-1}{k}{\kern-3pt{\textstyle} 2n}{}{k+1}\quad
-\qquad
\arbrebb{}{}{\raise4pt\hbox{$\scriptstyle
2n$}}{k}{\hbox{$\rond\atop k+1$}\hskip-5pt}{}{}{\raise5pt\hbox{$\scriptstyle\bullet\ k-1$}}
\cr
&:=E-F\cr
&=(E(\star)-F(\star))+(E(\rondsec)-F(\rondsec))\cr
\noalign{\bigskip\bigskip\smallskip}
&=F(\rondsec)=\arbrebbj{}{k+1}{}%
{\bullet}{\textstyle\rondsec}{2n}{k-1}{k}
\qquad
=\qquad\vbox{\vskip1.cm\offinterlineskip
\segment(0,0)\dir(0,1)\long{3}
\segment(0,3)\dir(-1,1)\long{4}
\segment(0,3)\dir(1,1)\long{12}
\segment(4,7)\dir(-1,1)\long{4}
\segment(8,11)\dir(-1,1)\long{4}
\centerput(-4,8.5){$\scriptstyle k-1$}
\centerput(0,11.5){$\scriptstyle 2n$}
\centerput(-4.7,6.5){$\bullet$}
\centerput(12,16){$\textstyle\rondsec$}
\centerput(6,6.5){$\scriptstyle k$}
\centerput(12,10.5){$\scriptstyle k+1$}
\centerput(4.5,15.5){$\carre$}
}\qquad\qquad.&(6.7)
\cr}
$$
By comparing the evaluations (6.6) and (6.7) we get:
$(B_{5}-D_{5})-(C_{3}-A_{5}-B_{2})=0$. This completes the proof of (6.5).\qed

\medskip
The entries $f_{2n}(3,2)$ and $f_{2n}(2n-1,2n-2)$, belonging to the upper diagonal of the lower triangle matrix~$M_{2n}$,
being evaluated by identity (6.4), and the entries of the upper triangle being known by Theorem~1.4, we can apply the crossing equalities (6.5), starting with $k=3$, and obtain all the values of the entries of that upper diagonal. Altogether, besides the entries of the upper triangle, all the entries lying on the border of the lower triangle can be calculated, as illustrated in boldface in the following matrix~$M_{8}$.

$$
M_{8}=\vcenter{\offinterlineskip
\halign{\vrule height10pt width.4pt
depth4pt$\hfil\ #\ $&\vrule height10pt width.4pt
depth4pt$\hfil\ #\ $&
$\hfil\ #\ $&$\hfil\ #\ $&$\hfil\ #\ $&$\hfil\ #\ $%
&$\hfil\ #\ $&$\hfil\ #\ $
&\vrule height10pt width.4pt depth4pt$\hfil\,#\,\hfil$%
\vrule height10pt width.4pt depth4pt\cr 
\noalign{\hrule}
k=&1&2&3&4&5&6&7&f_8(m,.)\cr 
\noalign{\hrule}
m=2 &       .  &   .  &  \bf 5  &  \bf  15 &  \bf   21 &  \bf   15 &    \bf  5 & \bf  61  \cr
 3&     \bf  5  &  \bf  \bf  \bf  10  &   . &  \bf   45 &    \bf  \bf 63 &  \bf   45 &  \bf   15& \bf 183  \cr
 4&   \bf   15  &  35  &  \bf  50  &   . &  \bf  101 &  \bf   63  &  \bf  21 & \bf 285\cr
 5&    \bf   21  &  54  &  86  & \bf  106  &   . &    \bf 45  &   \bf 15 & \bf 327   \cr
  6&     \bf   15 &   46 &   82 &   87 &   \bf  50 &    .  &   \bf  5 & \bf 285   \cr
7&    \bf  5  &  22  &  46  &  60 &   40 &   \bf  10  &   .  & \bf 183  \cr
8&      . &    \bf 16  &  \bf  16 &  \bf   14  &  \bf  10   &  \bf  5  &   . & \bf 61   \cr
\noalign{\hrule}
f_8(.,k)& \bf 61& \bf 183& \bf 285& \bf 327& \bf 285& \bf 183& \bf 61& \bf E_8=1385\cr
\noalign{\hrule}
}}$$

\smallskip
\centerline{Table 6.1: the matrix $M_{8}$, bold-faced entries analytically evaluated.}

\goodbreak
\bigskip
\centerline{\bf 7. Generating functions for the $f_{2n}(m,k)$}

\medskip
The calculation of the generating functions for the $f_{2n}(m,k)$
is similar to the calculation made for the tangent tree case in a our previous paper [FH13]. Recall the definition and some basic properties of the Poupard matrix. Let $G=(g_{i,j})$ $(i\ge 0,\,j\ge 0)$ be an infinite matrix with nonnegative integral entries. Say that~$G$ is a {\it Poupard matrix}, if for every $i\ge 0$, $j\ge 0$ the following identity holds:
$$
g_{i,j+2}-2\,g_{i+1,j+1}+g_{i+2,j}+4\,g_{i,j}=0.\leqno(7.1)
$$

{\it Remark}. Last coefficient is $4$, not $2$ as in [FH13].
\medskip

Let $G(x,y):=\sum\limits_{i\ge 0,\,j\ge 0}g_{i,j}\,(x^i/i!)\,(y^j/j!)$;\quad $R_{i}(y):=\sum\limits_{j\ge 0}g_{i,j}\,(y^j/j!)$ $(i\ge 0)$;\quad $C_{j}(x):=\sum\limits_{i\ge 0}g_{i,j}\,(x^i/i!)$ $(j\ge 0)$ be the exponential generating functions for the matrix itself, its rows and columns, respectively.
Propositions 9.1 and 9.2 in [FH13] can be rewritten as follows.

\proclaim Proposition 7.1. The following four properties are equivalent.\hfil\break\indent
(i) $G=(g_{i,j})$ $(i\ge 0,\,j\ge 0)$ is a Poupard matrix;
\hfil\break\indent
(ii) $R''_i(y)-2\,R'_{i+1}(y)+R_{i+2}(y) +4\,R_{i}(y)=0$
for all $i\ge 0$;
\hfil\break\indent
(iii) $C''_{j}(x)-2\,C'_{j+1}(x)+C_{j+2}(x)+4\,C_{j}(x)$ for all $j\ge 0$;
\hfil\break\indent
(iv) $G(x,y)$ satisfies the partial differential equation:
$$
{\partial^2 G(x,y)\over \partial x^2}
-2\,{\partial^2 G(x,y)\over \partial x\,\partial y}
+{\partial^2 G(x,y)\over \partial y^2}
+4\,G(x,y)=0.\leqno(7.2)
$$

\proclaim Proposition 7.2. Let $G(x,y)$ be the exponential generating function for a Poupard matrix~$G$. Then,  
$$
G(x,y)=A(x+y)\,\cos( 2\,y)+B(x+y)\,
\sin( 2\,y),\leqno(7.3)
$$
where $A(y)$ and $B(y)$ are two arbitrary series.

\medskip

\goodbreak
The entries of the {\it upper triangles} of the matrices $(M_{2n})$ (see Table 1.1) are now recorded as entries of infinite matrices
$(\Omega^{(p)})$ $(p\ge 1)$ by 
$$
\omega^{(p)}_{i,j}:=\cases{0,&if $i+j\equiv p$ mod 2;\cr
f_{2n}(m,k),&if $i+j\not\equiv p$ mod 2;\cr}\leqno(7.4)
$$
with $m:=p+1$, $k:=p+j+2$, $2n:=p+i+j+3$. Conversely, $i:=2n-k-1$, $j:=k-m-1$, $p:=m-1$. In particular,
the first one
$\Omega^{(1)}=(\omega_{i,j}^{(1)})$ $(i,j\ge 0)$ contains the first rows of the upper triangles, displayed as counter-diagonals. Furthermore, a counter-diagonal with zero entries is placed between two successive rows. 

\medskip
$\Omega^{(1)}={}$

{\eightpoint
\vskip-10pt
$$\eqalignno{
&\bordermatrix{&0&1&2&3&4&5&6&7\cr
0&f_{4}(2,3)&0&f_{6}(2,5)&0&f_{8}(2,7)&0&f_{10}(2,9)&0\cr
1&0&f_{6}(2,4)&0&f_{8}(2,6)&0&f_{10}(2,8)&0&\cdots\cr
2&f_{6}(2,3)&0&f_{8}(2,5)&0&f_{10}(2,7)&0&\cdots\cr
3&0&f_{8}(2,4)&0&f_{10}(2,6)&0&\cdots\cr
4&f_{8}(2,3)&0&f_{10}(2,5)&0&\cdots\cr
5&0&f_{10}(2,4)&0&\cdots\cr
6&f_{10}(2,3)&0&\cdots\cr
7&0&\cdots\cr
}\cr
&\ =\bordermatrix{&0&1&2&3&4&5&6&7\cr
0&1&0&1&0&5&0&61&0\cr
1&0&3&0&15&0&183&0&\cdots\cr
2&1&0&21&0&285&0&\cdots\cr
3&0&15&0&327&0&\cdots\cr
4&5&0&285&0&\cdots\cr
5&0&183&0&\cdots\cr
6&61&0&\cdots\cr
7&0&\cdots\cr
}.\cr
}
$$
}

\proclaim Proposition 7.3. Every matrix $\Omega^{(p)}$
$(p\ge 1)$ is a Poupard matrix.

\proof
Using Definition (7.4) we have
$$\displaylines{\quad
\omega^{(p)}_{i,j+2}-2\, \omega^{(p)}_{i+1,j+1}
+\omega^{(p)}_{i+2,j} +4\,\omega^{(p)}_{i,j}\hfill\cr
\kern2cm{}
=f_{2n+2}(m,k+2)-2\,f_{2n+2}(m,k+1)\hfill\cr
\kern5cm{}+f_{2n+2}(m,k)+4\,f_{2n}(m,k)
\hfill\cr
\kern2cm{}
=\Deltaa_{k}f_{2n+2}(m,k)+4\,f_{2n}(m,k)=0,\hfill\cr
}$$
by rule $(R\,2)$.\qed

The row labeled~$i$ of $\Omega^{(p)}$ will be denoted by 
$\Omega^{(p)}_{i, \brullet}$ and the exponential generating function for that 
row by
$\Omega^{(p)}_{i, \brullet}(y)=\sum_{j\geq 0} \omega_{i,j}^{(p)} {y^j/ j!}$. 
Also,
$\Omega^{(p)}(x,y):=\sum_{i\ge 0}\Omega^{(p)}_{i,\brullet}(y)x^i/i!$ will be the double exponential generating function for the matrix $\Omega^{(p)}$.
\proclaim Proposition 7.4. For all $p\ge 1$ we have:
$$
\Omega^{(p+1)}_{0,\brullet}(y)
=\Omega^{(1)}_{p,\brullet}(y)\hbox{\quad and \quad}
\Omega^{(p)}_{1,\brullet}(y)
=3\,{d\over dy}\Omega^{(p)}_{0,\brullet}(y)
=3\,{d\over dy}\Omega^{(1)}_{p-1,\brullet}(y).
$$

\proof
For the first identity it suffices to prove 
$\omega^{(p+1)}_{0,j}= \omega^{(1)}_{p,j}$, that is
$$
f_{2n}(m, 2n-1) = f_{2n}(2, 2n+1-m).
$$
This is true by the symmetry property of the Poupard triangle, as proved in Corollary 1.3 in [FH12].
For the second identity it suffices to prove
$\omega^{(p)}_{1,j}= 3 \omega^{(p)}_{0,j+1}$, that is
$$
f_{2n}(m, 2n-2) = 3f_{2n}(m, 2n-1).
$$
This is true by Theorem 1.2.\qed

\medskip
As $x$ and $y$ play a symmetric role in (7.3), the solution in 
(7.3) may also be written
$$\eqalignno{
G(x,y)&=A(x+y)\,\cos(2\,x)+B(x+y)\,\sin(2\,x),\cr
\noalign{\hbox{
so that the generating function of each matrix 
$\Omega^{(p)}$ is of the form}}
\Omega^{(p)}(x,y)
&=A(x+y)\,\cos( 2\,x)+B(x+y)\,\sin(2\,x).\cr}
$$
Hence, 
$\Omega^{(p)}(x,y)\Bigm|_{\{x=0\}}=\Omega^{(p)}_{0,\brullet}(y)=A(y)$. Also, 
$$
\leqalignno{({\partial/\partial x})\Omega^{(p)}(x,y)
	&=(({\partial/\partial x})A(x+y))\,\cos( 2\,x)
+A(x+y)(-2)\sin(2\,x)\cr
&\quad +(({\partial/\partial x})B(x+y))\,\sin( 2\,x)
+B(x+y)\,2\,\cos( 2\,x)\cr
}
$$ 
and
$$\leqalignno{
{\partial\over\partial x}\Omega^{(p)}(x,y)
\Bigm|_{\textstyle\{x=0\}}&=
{\partial\over\partial x}A(x+y)\Bigm|_{\textstyle\{x=0\}}
+2\,B(x+y)\Bigm|_{\textstyle\{x=0\}}\cr
&={d\over dy}A(y)
+2\,B(y)\cr
&=\Omega^{(p)}_{1,\brullet}(y).\cr}
$$
By Proposition 7.4 we have
$A(y)=\Omega^{(1)}_{p-1,\brullet}(y)$ and
$B(y)= (\Omega^{(p)}_{1,\brullet}(y) - \displaystyle{ d\over dy}A(y))/2
= \Bigl(3  \displaystyle{ d\over dy} \Omega^{(1)}_{p-1,\brullet}(y)- \displaystyle{ d\over dy}  \Omega^{(1)}_{p-1,\brullet}(y) \Bigr)/2
= \displaystyle{ d\over dy}  \Omega^{(1)}_{p-1,\brullet}(y)$.  Hence,
$$
\Omega^{(p)}(x,y)
=\cos(2\,x)\,\Omega^{(1)}_{p-1,\brullet}(x+y)
+\sin(2\,x)\,{\partial \over \partial y}\Omega^{(1)}_{p-1,\brullet}(x+y).\leqno(7.5)
$$ 

First, make the evaluation of $\Omega^{(1)}(x,y)$. The row labeled~0 of the matrix $\Omega^{(1)}$ reads:
1, 0, 1, 0, 5, 0, 61, \dots, which is the sequence of the coefficients of the Taylor  expansion of $\sec y$. Thus, 
$\Omega^{(1)}_{0,\brullet}(y)=\sec y$. Taking $p=1$ in (7.5) we get
$$
\leqalignno{
	\Omega^{(1)}(x,y)&=\sec(x+y)\cos(2\,x) +\sec (x+y)\,\tan (x+y)\,\sin(2\,x)\cr
	&={\cos (x-y)\over \cos^2(x+y)}. &(7.6)\cr
}
$$

\goodbreak
For further use let us also calculate the partial derivative of
$\Omega^{(1)}(x,y)$ with respect of~$y$:
$$
\leqalignno{
{\partial\over \partial y}
\Omega^{(1)}(x,y)
&={1\over \cos^3(x+y)}
\Bigl(\sin(x-y)\cos(x+y) +2\cos(x-y)\sin(x+y)\Bigr)\cr
&={1\over 2\,\cos^3(x+y)}\bigl(\sin(2y)+3\,\sin(2x)\bigr).
&(7.7)\cr
}
$$

Now, define
$$\Omega(x,y,z):=\displaystyle\sum_{p\ge 1}
\Omega^{(p)}(x,y){z^{p-1}\over (p-1)!}\leqno(7.8)
$$
and make use of (7.5)---(7.8):
$$
\leqalignno{
\Omega(x,y,z)
&=\cos(2x)\sum_{p\ge 1}
\Omega^{(1)}_{p-1,\brullet}(x+y){z^{p-1}\over (p-1)!}\cr
&\kern2cm{}
+\sin(2x)\,{\partial\over \partial y}
\Bigl(\sum_{p\ge 1}
\Omega^{(1)}_{p-1,\brullet}(x+y){z^{p-1}\over (p-1)!}\Bigr)\cr
&=\cos(2x)\,\Omega^{(1)}(z,x+y)+
\sin(2x)\,{\partial\over \partial y}\,\Omega^{(1)}(z,x+y)\cr
&=\cos(2x)\,{\cos(x+y-z)\over\cos^2(x+y+z)}\cr
&\kern2cm{}+\sin(2x)\,
{\sin(2(x+y))+3\,\sin(2z)\over 2\,\cos^3(x+y+z)}\cr
&={1\over 2\,\cos^3(x+y+z)}
\Bigl(\cos(2x)\,\bigl(\cos(2(x+y))+\cos(2z)\bigr)\cr
&\kern3.5cm{}+\sin(2x)\,\bigl(\sin(2(x+y))+3\sin(2z)\bigr)\Bigr)\cr
&={1\over 2\,\cos^3(x+y+z)}
\Bigl(\cos(2y)+2\,\cos(2(x-z))-\cos(2(z+x))\Bigr).&(7.9)\cr
}
$$

By definition of the 
$\omega^{(p)}_{i,j}$'s we get
$$
\leqalignno{
\Omega(x,y,z)&=\sum_{p,i,j}\omega^{(p)}_{i,j}
{z^{p-1}\over (p-1)!}
{x^{i}\over i!}
{y^{j}\over j!}\quad(p\ge 1,\,i\ge 0,\,j\ge 0);\cr
&=\sum_{k,m,n}
f_{2n}(m,k){x^{2n-k-1}\over (2n-k-1)!}
{y^{k-m-1}\over (k-m-1)!}
{z^{m-2}\over (m-2)!},&(7.10)\cr
 }
$$
the latter sum over the set
$\{3\le m+1\le k\le 2n-1\}$. We have proved Theorem~1.5.

\vfill\eject
\vglue2cm

\centerline{\bf References}

{\eightpoint

\bigskip

\article  An1879|D\'esir\'e Andr\'e|D\'eveloppement de $\sec x$ et
$\tan x$|C. R. Math. Acad. Sci. Paris|88|1879|965--979|

\article An1881|D\'esir\'e Andr\'e|Sur les permutations
altern\'ees|J. Math. Pures et Appl.|7|1881|167--184|

\livre Co74|Louis Comtet|\hskip-5pt Advanced 
Combinatorics|\hskip-5pt D.
Reidel/Dordrecht-Holland, Boston, {\oldstyle 1974}|

\article En66|R.C. Entringer|A combinatorial interpretation of the
Euler and Bernoulli numbers|Nieuw Arch. Wisk.|14|1966|241--246|

\divers FH12|Dominique Foata; Guo-Niu Han|Finite Difference Calculus
for Alternating Permutations, preprint, 15~p|

\divers FH13|Dominique Foata; Guo-Niu Han|Tree Calculus for Bivariable Difference Equations, preprint, 36~p|

\article GHZ10|Yoann Gelineau; Heesung Shin; Jiang Zeng|Bijections for
Entringer families|Europ. J. Combin.|32|2011|100--115|

\divers Ha12|Guo-Niu Han|The Poupard Statistics on Tangent and Secant Trees, preprint 12~p|

\livre Jo39|Charles Jordan|Calculus of Finite Differences|R\"ottig and
Romwalter,  Budapest, {\oldstyle 1939}|

\article KPP94|A. G. Kuznetsov; I. M. Pak;  A. E.
Postnikov|Increasing trees and alternating permutations|Uspekhi Mat.
Nauk|49|1994|79--110|

\livre Ni23|Niels Nielsen|Trait\'e \'el\'ementaire des nombres
de Bernoulli|Paris, Gauthier-Villars, {\oldstyle 1923}|

\article Po82|Christiane Poupard|De nouvelles significations
\'enum\'eratives des nombres\hfil\break d'Entringer|Discrete
Math.|38|1982|265--271|

\article Po89|Christiane Poupard|Deux propri\'et\'es des arbres
binaires ordonn\'es stricts|Europ. J. Combin.|10|1989|369--374|

\article  Po97|Christiane Poupard|Two other interpretations of the
Entringer numbers|Europ. J. Combin.|18|1997|939--943|

\divers Sl07|N.J.A. Sloane|On-line Encyclopedia of Integer
Sequences,\hfil\break
{\tt
http://www.research.att.com/\char126 njass/sequences/}|

\divers St10|Richard P.
Stanley|A Survey of Alternating Permutations, in {\sl Combinatorics and graphs}, 165--196, {\sl Contemp. Math.}, {\bf 531}, Amer. Math. Soc. Providence, RI, {\oldstyle 2010}|

\divers Vi88|Xavier G. Viennot|S\'eries g\'en\'eratrices
\'enum\'eratives, chap.~3, Lecture Notes, 160~p., 1988, notes de
cours donn\'es
\`a l'\smash{\'E}cole Normale Sup\'erieure Ulm (Paris), UQAM (Montr\'eal,
Qu\'ebec) et Universit\'e de Wuhan (Chine)\hfil\break
{\tt
http://web.mac.com/xgviennot/Xavier\_Viennot/cours.html}|

\bigskip\bigskip
\hbox{\vtop{\halign{#\hfil\cr
Dominique Foata \cr
Institut Lothaire\cr
1, rue Murner\cr
F-67000 Strasbourg, France\cr
\noalign{\smallskip}
{\tt foata@unistra.fr}\cr}}
\qquad
\vtop{\halign{#\hfil\cr
Guo-Niu Han\cr
I.R.M.A. UMR 7501\cr
Universit\'e de Strasbourg et CNRS\cr
7, rue Ren\'e-Descartes\cr
F-67084 Strasbourg, France\cr
\noalign{\smallskip}
{\tt guoniu.han@unistra.fr}\cr}}}

}

\bye